arXiv:1707.06551v2 [math.NA] 9 Feb 2018

# Boundary integral equation analysis for suspension of spheres in Stokes flow

Eduardo Corona * and Shravan Veerapaneni *

**Abstract**

We show that the standard boundary integral operators, defined on the unit sphere, for the Stokes equations diagonalize on a specific set of vector spherical harmonics and provide formulas for their spectra. We also derive analytical expressions for evaluating the operators away from the boundary. When two particle are located close to each other, we use a truncated series expansion to compute the hydrodynamic interaction. On the other hand, we use the standard spectrally accurate quadrature scheme to evaluate smooth integrals on the far-field, and accelerate the resulting discrete sums using the fast multipole method (FMM). We employ this discretization scheme to analyze several boundary integral formulations of interest including those arising in porous media flow, active matter and magneto-hydrodynamics of rigid particles. We provide numerical results verifying the accuracy and scaling of their evaluation.

## 1 Introduction

Suspension of spherical particles in Stokes flow acts as a mimetic model for several natural and engineering systems. Often, the physical phenomena of interest happen at scales much larger than the constituent particle sizes e.g., collective motion in bacterial suspensions [10], bulk rheology of polydisperse colloidal suspensions, pattern formations in electro- and magneto-rheological fluids [22] and self-assembly of particles [27]. Consequently, direct simulation methods that scale to large number of particles and that are numerically stable for long-time simulations are crucial to enable insights into these systems.

Several techniques have been developed in the past few decades for simulating the hydrodynamics of multiple spherical particles including the Stokesian dynamics approach (e.g., [4, 11, 47]), multipole methods (e.g., [8, 37]), fictitious domain methods (e.g., [33]) and boundary integral methods (e.g., [1, 9, 34]). We refer the reader to [29] for a recent review on the broader topic of simulation methods for particulate flows. The present work combines features from both the multipole methods (spectral representations) and the boundary integral methods (second-kind formulations, fast algorithms) to arrive at a fast, spectrally accurate numerical method.

Our work is closely related to three recent efforts, that of Veerapaneni et al. [41], Vico et al. [43] and Singh et al. [38] (listed in chronological order). In [41], using the antenna theorems of [37], the spectrum of the "single-layer" Stokes boundary integral operator (BIO) was derived and applied to analyze certain integro-dfferential operators on the sphere; here, we extend this framework to all the other relevant BIOs. In [43], authors derived signatures of the BIOs that arise when solving Helmholtz or Maxwell equations in the frequency domain. The present work can be viewed as an extension of [43] to the Stokes equations. In [38], a Galerkin approach for evaluating Stokes BIOs on spheres was developed. The main difference from the present work is the choice of the basis functions: while *tensorial spherical harmonics* were used in [38], we chose a specific set of *vector spherical harmonics*. We show that this choice leads to much simpler formulas and diagonalization of most of the BIOs. Another important distinction is that the $N$-body hydrodynamic interactions of the particles were computing directly in [38], leading to a quadratic complexity in $N$ whereas our method is linear in $N$ via the use of Nyström's

*Department of Mathematics, University of Michigan

method for evaluating the smooth far-field integrals, accelerated by the fast multipole method (FMM) [16, 18, 39].

*Synopsis.* We consider several BIOs that arise when solving the Stokes equations and compute their signatures analytically on the unit sphere. This enables us to convert the classical task of *weakly-singular* integral evaluation to simple formula evaluation. We also present formulas for evaluating the operators at arbitrary target locations away from the sphere. Thereby, the issue of accurate evaluation of *nearly-singular* integrals also reduces to simple analytic expression evaluation. We then demonstrate the solution procedure for various physical problems using the standard integral equation formulations proposed in the literature.

The paper is organized as follows. In §2, we introduce the basis functions for representing scalar and vector fields on the sphere and the definitions of the boundary integral operators. In §3, we derive the signatures of these operators and the analytical formulas for evaluating the velocity and pressure away from a unit sphere given a certain form of the jump conditions. We use these formulas in §4 to develop a fast, spectrally accurate singular and nearly-singular integral evaluation scheme. Finally, in §5, we discuss the standard model problems in creeping flow of spherical particle suspensions, their reformulation as boundary integral equations and perform a series of numerical experiments to validate our solvers.

# 2 Mathematical Preliminaries

In this section, we provide a summary of the spherical harmonic basis functions, which will be used to represent scalar and vector fields defined on the sphere, and also provide definitions for the classical boundary integral operators that arise when solving the Laplace and Stokes equations. The concepts discussed here are fairly standard e.g., see [21, 32, 35, 42].

## 2.1 Spherical harmonic bases

### Scalar spherical harmonics

**Definition 2.1.** *Let $\theta$ and $\phi$ be the polar and azimuthal angles in the standard parametrization of the unit sphere. The scalar spherical harmonic $Y_n^m$ of degree $n$ and order $m$ (for $|m| \leq n$) is defined in terms of the associated Legendre functions $P_n^m$ by*

$$Y_n^m(\theta, \phi) = \sqrt{\frac{2n+1}{4\pi}} \sqrt{\frac{(n-|m|)!}{(n+|m|)!}} \, P_n^{|m|}(\cos\theta) \, e^{im\phi}. \tag{2.1}$$

Scalar spherical harmonics form an orthonormal basis of eigenfunctions of the Laplacian for square-integrable functions on the unit sphere. That is, any function $\sigma \in L^2(\mathbb{S}^2)$ has the expansion:

$$\sigma(\theta, \phi) = \sum_{n=0}^{\infty} \sum_{m=-n}^{n} \widehat{\sigma}_n^m \, Y_n^m(\theta, \phi), \quad \theta \in [0, \pi], \quad \phi \in [0, 2\pi], \tag{2.2}$$

$$\text{where} \quad \widehat{\sigma}_n^m = \int_0^{2\pi} \int_0^{\pi} \sigma(\theta, \phi) \overline{Y_n^m(\theta, \phi)} \, \sin\theta \, d\theta d\phi. \tag{2.3}$$

For $\sigma \in C^\infty(\mathbb{S}^2)$, the finite-term approximation (truncating the outer sum to $n = 0 \ldots p$ in Eq. (2.2), which yields $(p+1)^2$ terms) is spectrally convergent with $p$ [32]. One can use fast transforms for both the longitude (Fast Fourier Transform or FFT) and the latitude (Fast Legendre Transform or FLT) to implement a fast forward spherical harmonic transform that computes the coefficients for the approximation of order $p$ in $\mathcal{O}(p^2 \log^2 p)$ operations. The inverse transform can be obtained in a similar fashion [30]. However, note that the break-even point for existing FLTs is large and typically, only the FFTs will be employed with a complexity of $\mathcal{O}(p^3 \log p)$ for forward and inverse transforms.

### Vector spherical harmonics

Vector spherical harmonics are an extension of the scalar spherical harmonics to square-integrable vector fields on the sphere. They can in fact be defined in terms of scalar spherical harmonics and their derivatives.

**Definition 2.2.** *The vector spherical harmonics* $\boldsymbol{V}_n^m, \boldsymbol{W}_n^m$ *and* $\boldsymbol{X}_n^m$ *of degree* $n$ *and order* $m$ *(for* $|m| \le n$*) are defined by*

$$\boldsymbol{V}_n^m = \nabla_\gamma Y_n^m(\theta,\phi) - (n+1) Y_n^m(\theta,\phi)\boldsymbol{e}_r(\theta,\phi), \tag{2.4}$$

$$\boldsymbol{W}_n^m = \nabla_\gamma Y_n^m(\theta,\phi) + n Y_n^m(\theta,\phi)\boldsymbol{e}_r(\theta,\phi), \tag{2.5}$$

$$\boldsymbol{X}_n^m = \boldsymbol{e}_r(\theta,\phi) \times \nabla_\gamma Y_n^m(\theta,\phi), \tag{2.6}$$

*where* $\nabla_\gamma = \frac{1}{\sin\theta}\frac{\partial}{\partial\phi}\boldsymbol{e}_\phi + \frac{\partial}{\partial\theta}\boldsymbol{e}_\theta$ *is the surface gradient operator, and* $\boldsymbol{e}_r, \boldsymbol{e}_\theta, \boldsymbol{e}_\phi$ *are the radial, polar and azimuthal unit vectors. For the sake of notational convenience, we will suppress the dependency on* $(\theta,\phi)$*.*

Depending on the problem or operator of interest, various sets of vector spherical harmonics have been employed in the literature. The ones defined above were first proposed in [21]. One notable alternative [2, 31] may be obtained from the split of vector fields $\boldsymbol{F} \in L^2(\mathbb{S}^2)^3$ into a radial component $f_r\boldsymbol{e}_r$ and a tangential component $f_\theta\boldsymbol{e}_\theta + f_\phi\boldsymbol{e}_\phi$. Radial vector fields are spanned by the orthogonal basis $\boldsymbol{Y}_n^m = Y_n^m\boldsymbol{e}_r$, while tangential fields are spanned by the orthogonal basis consisting of $\boldsymbol{G}_n^m = \nabla_\gamma Y_n^m$ and $\boldsymbol{X}_n^m$ (defined in (2.6)), for $n \ge 0, |m| \le n$.

### Vector spherical harmonic transforms

As suggested in [24], we compute the vector spherical harmonic transform and its inverse by means of three fast scalar transforms, one for each coordinate vector $\boldsymbol{e}_r, \boldsymbol{e}_\theta, \boldsymbol{e}_\phi$. Given a vector field $\boldsymbol{F} = f_r\boldsymbol{e}_r + f_\theta\boldsymbol{e}_\theta + f_\phi\boldsymbol{e}_\phi$, an expansion of $f_r$ in scalar spherical harmonics trivially corresponds to an expansion of the radial part in $Y_n^m\boldsymbol{e}_r$, with coefficients $\varphi_{n,m}^Y$.

In order to obtain the corresponding scalar coefficients for the tangential part, which we denote $\{\varphi_{n,m}^\theta, \varphi_{n,m}^\phi\}$, we note that using standard recurrence relations for Legendre functions, we can derive the following relations between the bases $\{\boldsymbol{G}_n^m, \boldsymbol{X}_n^m\}$ and $\{Y_n^m\boldsymbol{e}_\theta, Y_n^m\boldsymbol{e}_\phi\}$:

$$n_1(n+1)(\sin\theta)Y_n^m\boldsymbol{e}_\theta = \left(-\alpha_n^m\boldsymbol{G}_{n+1}^m + \beta_n^m\boldsymbol{G}_{n-1}^m\right) - im\boldsymbol{X}_n^m, \tag{2.7}$$

$$n_1(n+1)(\sin\theta)Y_n^m\boldsymbol{e}_\phi = \left(-\alpha_n^m\boldsymbol{X}_{n+1}^m + \beta_n^m\boldsymbol{X}_{n-1}^m\right) + im\boldsymbol{G}_n^m, \tag{2.8}$$

as well as the reverse (for $n > 0$):

$$(\sin\theta)\boldsymbol{G}_n^m = \left(\alpha_n^m Y_{n+1}^m - \beta_n^m Y_{n-1}^m\right)\boldsymbol{e}_\theta + imY_n^m\boldsymbol{e}_\phi, \tag{2.9}$$

$$(\sin\theta)\boldsymbol{X}_n^m = \left(\alpha_n^m Y_{n+1}^m - \beta_n^m Y_{n-1}^m\right)\boldsymbol{e}_\phi - imY_n^m\boldsymbol{e}_\theta, \tag{2.10}$$

where $n_1 = \max(n,1)$, $\alpha_n^m = \sqrt{\frac{n^2((n+1)^2-m^2)}{(2n+1)(2n+3)}}$ and $\beta_n^m = \sqrt{\frac{(n+1)^2(n^2-m^2)}{(2n+1)(2n-1)}}$, and undefined terms (e.g., $Y_{-1}^m$ or $Y_n^{n+1}$) are dropped. Given coefficients $\{\varphi_{n,m}^\theta, \varphi_{n,m}^\phi\}$ from scalar spherical harmonic expansions of $f_\theta/\sin\theta$ and $f_\phi/\sin\theta$, using Eq. (2.7) and Eq. (2.8) we get

$$\begin{bmatrix}\varphi_{n,m}^G \\ \varphi_{n,m}^X\end{bmatrix} = T_{\text{fwd}}\begin{bmatrix}\varphi_{n,m}^\theta \\ \varphi_{n,m}^\phi\end{bmatrix}, \tag{2.11}$$

with $T_{\text{fwd}}$ being sparse and permutable into a pentadiagonal matrix. Finally, given an expansion of a vector field in the vector spherical harmonic basis $\{\boldsymbol{Y}_n^m, \boldsymbol{G}_n^m, \boldsymbol{X}_n^m\}$,

$$\boldsymbol{F} = \sum_{n,m} \varphi_{n,m}^Y\boldsymbol{Y}_n^m + \varphi_{n,m}^G\boldsymbol{G}_n^m + \varphi_{n,m}^X\boldsymbol{X}_n^m, \tag{2.12}$$

by the definitions of $\boldsymbol{V}_n^m$ and $\boldsymbol{W}_n^m$ (equations (2.4) and (2.5)), we can find the corresponding coefficients $\varphi^V_{n,m}, \varphi^W_{n,m}$ for the expansion of $\boldsymbol{F}$ in our original set of vector spherical harmonics, with

$$\varphi^V_{n,m} = (n\varphi^G_{n,m} - \varphi^Y_{n,m})/(2n+1), \tag{2.13}$$

$$\varphi^W_{n,m} = ((n+1)\varphi^G_{n,m} + \varphi^Y_{n,m})/(2n+1). \tag{2.14}$$

**Inverse transforms.** Given coefficients of a vector field in any of our two vector spherical harmonic bases, it is possible to relate them to scalar expansions in the radial and tangential directions. The algorithm to evaluate the corresponding vector field thus proceeds in two stages: conversion to scalar coefficients $\{\varphi^Y_{n,m}, \varphi^\theta_{n,m}, \varphi^\phi_{n,m}\}$ and the application of three scalar inverse transforms to evaluate the field's spherical coordinate components.

Using Eq. (2.9) and Eq. (2.10), we similarly derive the "inverse" of Eq. (2.11):

$$\begin{bmatrix} \varphi^\theta_{n,m} \\ \varphi^\phi_{n,m} \end{bmatrix} = T_{\text{inv}} \begin{bmatrix} \varphi^G_{n,m} \\ \varphi^X_{n,m} \end{bmatrix}. \tag{2.15}$$

If we wish to convert from $\{\boldsymbol{V}_n^m, \boldsymbol{W}_n^m, \boldsymbol{X}_n^m\}$, we first invert Eq. (2.13) and Eq. (2.14) to get:

$$\varphi^G_{n,m} = \varphi^V_{n,m} + \varphi^W_{n,m}, \tag{2.16}$$

$$\varphi^Y_{n,m} = -(n+1)\varphi^V_{n,m} + n\varphi^W_{n,m}. \tag{2.17}$$

## 2.2 Boundary integral operators

Let $\boldsymbol{n}(\boldsymbol{y})$ be the unit outward normal vector at $\boldsymbol{y} \in \Gamma = \mathbb{S}^2$. The fundamental solution to the Laplace equation in free space, denoted by $G^L$ is given by, $G^L(\boldsymbol{x}, \boldsymbol{y}) = \frac{1}{4\pi|\boldsymbol{x}-\boldsymbol{y}|}$. The Laplace single and double layer potentials are then defined respectively as

$$\mathcal{S}^L[\sigma](\boldsymbol{x}) = \int_\Gamma G^L(\boldsymbol{x}, \boldsymbol{y})\sigma(\mathbf{y})d\Gamma(\mathbf{y}) \quad \text{and} \quad \mathcal{D}^L[\mu](\boldsymbol{x}) = \int_\Gamma \frac{\partial}{\partial \boldsymbol{n}(\boldsymbol{y})} G^L(\boldsymbol{x}, \boldsymbol{y})\mu(\mathbf{y})d\Gamma(\mathbf{y}), \tag{2.18}$$

where $\sigma, \mu$ in $L^2(\Gamma)$ are often termed as the corresponding "*density functions*" or simply as "*densities*" of the layer potentials. Operators on the surface $\Gamma$ are defined by taking limits as $\mathbf{x} \to \Gamma$ along the normal direction. For the single-layer, standard jump conditions show that interior and exterior limits coincide [25], yielding a weakly singular operator. For the double-layer, the interior and exterior limits, which we denote by $D^L_+$ and $D^L_-$, differ and the jump across the surface is

$$[[\mathcal{D}^L\mu]]_\Gamma = \mathcal{D}^L_+[\mu] - \mathcal{D}^L_-[\mu] = \mu.$$

Let $G_{i,j}(\boldsymbol{x}, \boldsymbol{y})$ be the Stokeslet, that is, the fundamental solution to the Stokes equations in free space in $\mathbb{R}^3$, given by

$$G_{i,j}(\boldsymbol{x}, \boldsymbol{y}) = \frac{1}{8\pi}\left(\frac{\delta_{i,j}}{|\boldsymbol{x}-\boldsymbol{y}|} + \frac{(x_i - y_i)(x_j - y_j)}{|\boldsymbol{x}-\boldsymbol{y}|^3}\right), \tag{2.19}$$

and let $T_{i,j,k}$ be the Stresslet, also known as the traction kernel, given by

$$T_{i,j,k}(\boldsymbol{x}, \boldsymbol{y}) = -\frac{3}{4\pi}\frac{(x_i - y_i)(x_j - y_j)(x_k - y_k)}{|\boldsymbol{x}-\boldsymbol{y}|^5}. \tag{2.20}$$

Let $\boldsymbol{\mu}, \boldsymbol{\sigma}$ be the density functions in $L^2(\Gamma)^3$. Then, the single and double layer potentials for the Stokes equation are defined, similar to the Laplace case, as

$$\mathcal{S}[\boldsymbol{\sigma}](\boldsymbol{x})_i = \int_\Gamma G_{i,j}(\boldsymbol{x}, \boldsymbol{y})\boldsymbol{\sigma}_j(\boldsymbol{y})\, d\Gamma(\mathbf{y}) \quad \text{and} \quad \mathcal{D}[\boldsymbol{\mu}](\boldsymbol{x})_i = \int_\Gamma T_{i,j,k}(\boldsymbol{x}, \boldsymbol{y})\boldsymbol{n}_k(\boldsymbol{y})\boldsymbol{\mu}_j(\boldsymbol{y})\, d\Gamma(\mathbf{y}), \tag{2.21}$$

which again are well-defined for $\mathbf{x} \notin \Gamma$. Finally, the traction associated to the single layer $\mathcal{S}[\boldsymbol{\sigma}]$ is given by

$$\mathcal{K}[\boldsymbol{\sigma}](\boldsymbol{x})_i = \int_\Gamma T_{i,j,k}(\boldsymbol{x}, \boldsymbol{y})\boldsymbol{n}_k(\boldsymbol{x})\boldsymbol{\sigma}_j(\boldsymbol{y})\, d\Gamma(\mathbf{y}). \tag{2.22}$$

Definitions of these operators on the surface $\Gamma$, as well as jump conditions are in a sense vector analogues of those for the Laplace operator [35].

# 3 Layer potential evaluation and spectra

Given densities $\sigma, \mu \in L^2(\mathbb{S}^2)$ on the sphere, the single and the double layer potentials for the Laplace and Stokes equations constitute solutions to the corresponding equation away from the sphere, and satisfy growth conditions at infinity automatically. Using a separation of variables argument, e.g. for a solution $\varphi$ to the Laplace equation on the exterior or interior of the sphere, it is also well known that one may find expansions of the form

$$\varphi = \sum_{n=0}^{\infty} \sum_{m=-n}^{n} a_n^m f_n(r) Y_n^m(\theta, \phi). \tag{3.1}$$

Using the appropriate jump conditions and orthogonality of the spherical harmonics basis, this allows us to conclude that all the layer potentials defined on the sphere above diagonalize on this basis. For each member of the appropriate spherical harmonic basis, we present formulas to evaluate the layer potentials defined in Eq. (2.18) and Eq. (2.21) at an arbitrary $\boldsymbol{x} \notin \mathbb{S}^2$ with spherical coordinates $(r, \theta, \phi)$. We note that, by superposition, given an integral density $\sigma$ with the expansion in spherical harmonics in Eq. (2.2), these formulas allow an expansion of $\varphi = S^L[\sigma]$ of the form

$$\mathcal{S}^L[\sigma](r, \theta, \phi) = \sum_{n=0}^{\infty} \sum_{m=-n}^{n} \lambda_n^m \hat{\sigma}_{n,m} f_n(r) Y_n^m(\theta, \phi) \tag{3.2}$$

where $\lambda_n^m$ are the eigenvalues for the operator $\mathcal{S}^L$ defined on the sphere.

The formulas and operator spectra presented below are well-known, and part of a vast literature of analysis of integral equation operators on the sphere. For discussion on the analysis of Laplace, Helmholtz and Maxwell layer potentials and related operators, see [43]. For the Stokes equations, the eigenvalues and eigenfunctions of the single-layer operator were derived in [37]. The eigenvalues of the double-layer operator were also derived in several works using tensorial spherical harmonics as basis functions (e.g., Chapter 17 of [23] and [44], [46]). Here, we show that a single set of vectorial spherical harmonics are the eigenfunctions of all the BIOs required in practical applications. We include an outline of this procedure in Appendix §A.

## 3.1 Laplace equation

For the Laplace equation, proposing solutions $u_n^m(r, \theta, \phi) = f_n(r) Y_n^m(\theta, \phi)$ produces a homogeneous ODE for $f_n(r)$ with two sets of admissible solutions, valid for the exterior and interior problems, respectively. We first consider $\varphi_n^m = \mathcal{S}^L[Y_n^m]$. Although $\varphi_n^m$ is continuous across the sphere, jump conditions indicate that $[[\frac{\partial}{\partial r}\varphi_n^m]]_\Gamma = -Y_n^m$. Taking inner products with $Y_n^m$ yields

$$\mathcal{S}^L[Y_n^m](r, \theta, \phi) = \begin{cases} \frac{1}{2n+1} Y_n^m(\theta, \phi) r^{-n-1} & r \geq 1, \\ \frac{1}{2n+1} Y_n^m(\theta, \phi) r^n & r \leq 1. \end{cases} \tag{3.3}$$

Applying the same procedure to $\psi_n^m = \mathcal{D}^L[Y_n^m]$, using the jump conditions $[[\psi_n^m]]_\Gamma = Y_n^m$, $[[\frac{\partial}{\partial r}\psi_n^m]]_\Gamma = 0$, we obtain

$$\mathcal{D}^L[Y_n^m](r, \theta, \phi) = \begin{cases} \frac{n}{2n+1} Y_n^m(\theta, \phi) r^{-n-1} & r > 1, \\ -\frac{n+1}{2n+1} Y_n^m(\theta, \phi) r^n & r < 1. \end{cases} \tag{3.4}$$

**Theorem 1** (Laplace operators spectra). *On the unit sphere, both Laplace single and double layer operators diagonalize in the scalar spherical harmonics basis $Y_n^m$ and the spectra are given by*

$$\mathcal{S}^L[Y_n^m] = \frac{1}{2n+1} Y_n^m, \quad \mathcal{D}_+^L[Y_n^m] = \frac{n}{2n+1} Y_n^m, \quad \text{and} \quad \mathcal{D}_-^L[Y_n^m] = -\frac{(n+1)}{2n+1} Y_n^m. \tag{3.5}$$

## 3.2 Stokes equations

For the Stokes equations, we can apply a similar procedure by proposing ansatz velocity fields of the form $\boldsymbol{u}_n^m = f_n(r)\boldsymbol{V}_n^m + g_n(r)\boldsymbol{W}_n^m + h_n(r)\boldsymbol{X}_n^m$, with associated pressures of the form $p_n^m = \frac{q_n(r)}{r}Y_n^m$, which satisfy the Stokes equations. Enforcing the velocity and traction jump conditions for the single layer potential, we can obtain the following expressions

$$\mathcal{S}[\boldsymbol{V}_n^m](r,\theta,\phi) = \begin{cases} \dfrac{n}{(2n+1)(2n+3)}\boldsymbol{V}_n^m r^{-n-2} & r \geq 1, \\ \dfrac{n}{(2n+1)(2n+3)}\boldsymbol{V}_n^m r^{n+1} + \dfrac{n+1}{4n+2}\boldsymbol{W}_n^m(r^{n-1} - r^{n+1}) & r \leq 1. \end{cases} \tag{3.6}$$

$$\mathcal{S}[\boldsymbol{W}_n^m](r,\theta,\phi) = \begin{cases} \dfrac{n+1}{(2n+1)(2n-1)}\boldsymbol{W}_n^m r^{-n} + \dfrac{n}{4n+2}\boldsymbol{V}_n^m(r^{-n-2} - r^{-n}) & r \geq 1, \\ \dfrac{n+1}{(2n+1)(2n-1)}\boldsymbol{W}_n^m r^{n-1} & r \leq 1. \end{cases} \tag{3.7}$$

$$\mathcal{S}[\boldsymbol{X}_n^m](r,\theta,\phi) = \begin{cases} \dfrac{1}{2n+1}\boldsymbol{X}_n^m r^{-n-1} & r \geq 1, \\ \dfrac{1}{2n+1}\boldsymbol{X}_n^m r^{n} & r \leq 1. \end{cases} \tag{3.8}$$

The associated pressures are zero except for two cases: $w^+ = \mathcal{S}[\boldsymbol{W}_n^m]$ corresponds to $p = nY_n^m r^{-n-1}$ for $r \geq 1$, and $v^- = \mathcal{S}[\boldsymbol{W}_n^m]$ corresponds to $p = (n+1)Y_n^m r^n$ for $r \leq 1$.

Applying the same procedure to the double-layer potentials $\mathcal{D}[\boldsymbol{V}_n^m], \mathcal{D}[\boldsymbol{W}_n^m]$ and $\mathcal{D}[\boldsymbol{X}_n^m]$, and enforcing velocity and traction jumps across the sphere,

$$\mathcal{D}[\boldsymbol{V}_n^m](r,\theta,\phi) = \begin{cases} \dfrac{2n^2+4n+3}{(2n+1)(2n+3)}\boldsymbol{V}_n^m r^{-n-2} & r > 1, \\ \dfrac{-2n(n+2)}{(2n+1)(2n+3)}\boldsymbol{V}_n^m r^{n+1} + \dfrac{(n+1)(n+2)}{2n+1}\boldsymbol{W}_n^m(r^{n+1} - r^{n-1}) & r < 1. \end{cases} \tag{3.9}$$

$$\mathcal{D}[\boldsymbol{W}_n^m](r,\theta,\phi) = \begin{cases} \dfrac{2(n+1)(n-1)}{(2n+1)(2n-1)}\boldsymbol{W}_n^m r^{-n} + \dfrac{2n(n-1)}{4n+2}\boldsymbol{V}_n^m(r^{-n-2} - r^{-n}) & r > 1, \\ \dfrac{-(2n^2+1)}{(2n+1)(2n-1)}\boldsymbol{W}_n^m r^{n-1} & r < 1. \end{cases} \tag{3.10}$$

$$\mathcal{D}[\boldsymbol{X}_n^m](r,\theta,\phi) = \begin{cases} \dfrac{n-1}{2n+1}\boldsymbol{X}_n^m r^{-n-1} & r > 1, \\ \dfrac{-(n+2)}{2n+1}\boldsymbol{X}_n^m r^{n} & r < 1. \end{cases} \tag{3.11}$$

Again, there are only two cases for which the associated pressure is non-zero: $w^+ = \mathcal{D}[\boldsymbol{W}_n^m]$ corresponds to $p = n(n+1)Y_n^m r^{-n-2}$ for $r \geq 1$ and $v^- = \mathcal{D}[\boldsymbol{V}_n^m]$ corresponds to $p = -n(n+1)Y_n^m r^{n-1}$ for $r \leq 1$.

**Theorem 2** (Stokes operators spectra)**.** *On the unit sphere, both Stokes single and double layer operators, as well as the single layer traction operator diagonalize in the vector spherical harmonics basis* $\{\boldsymbol{V}_n^m, \boldsymbol{W}_n^m, \boldsymbol{X}_n^m\}$*. The corresponding eigenvalues are*

| | $\boldsymbol{V}_n^m$ | $\boldsymbol{W}_n^m$ | $\boldsymbol{X}_n^m$ |
|---|---|---|---|
| $\mathcal{S}$ | $\frac{n}{(2n+1)(2n+3)}$ | $\frac{n+1}{(2n+1)(2n-1)}$ | $\frac{1}{2n+1}$ |
| $\mathcal{D}_+$ | $\frac{2n^2+4n+3}{(2n+1)(2n+3)}$ | $\frac{2(n-1)(n+1)}{(2n+1)(2n-1)}$ | $\frac{n-1}{2n+1}$ |
| $\mathcal{D}_-$ | $\frac{-2n(n+2)}{(2n+1)(2n+3)}$ | $\frac{-(2n^2+1)}{(2n+1)(2n-1)}$ | $\frac{-(n+2)}{2n+1}$ |
| $\mathcal{K}_+$ | $\frac{-2n(n+2)}{(2n+1)(2n+3)}$ | $\frac{-(2n^2+1)}{(2n+1)(2n-1)}$ | $\frac{-(n+2)}{2n+1}$ |
| $\mathcal{K}_-$ | $\frac{2n^2+4n+3}{(2n+1)(2n+3)}$ | $\frac{2(n-1)(n+1)}{(2n+1)(2n-1)}$ | $\frac{n-1}{2n+1}$ |

(3.12)

As is the case for the Laplace operator, the spectra of $\mathcal{K}_+$ matches that of $\mathcal{D}_-$ and $\mathcal{K}_-$ that of $\mathcal{D}_+$. However, the fact that the normal vector in Eq. (2.22) depends on the target point $\boldsymbol{x}$ may introduce additional terms when attempting to evaluate this operator off the surface, as compared to the formulas in Eqs 3.6 - 3.11. In §3.3, we outline a general evaluation procedure for the evaluation of derivatives of BIOs of interest, and present the traction of a flow expanded in vector spherical harmonics as a relevant example.

## 3.3 Evaluation of derivatives

Once we define a solution to the Laplace or Stokes equations as a combination of single and double layer potentials, a boundary integral equation must be solved in order to find integral densities that match the given boundary conditions. For a suspension of multiple bodies, the formulas presented above allow us to map an expansion in spherical harmonics for an integral density $\sigma$ to the corresponding expansions for $\mathcal{S}[\sigma]$ and $\mathcal{D}[\sigma]$, and evaluate them on and off the surface of the unit sphere. We extend this to two examples of interest: the normal derivative or flux of a Laplace layer potential, and the traction force of a Stokes layer potential. We note that the procedure outlined below to numerically compute coefficients for evaluation of derivatives of layer potentials on and off the surface may be generalized to any sum of differential operators separable into radial and tangential parts in spherical coordinates.

**Flux calculation (Laplace).** Given $\varphi$ a linear combination of $\mathcal{S}^L[\sigma]$ and $\mathcal{D}^L[\sigma]$, Eqs. (3.3 - 3.4) provide a formula of the form,

$$\varphi(\boldsymbol{x}) = \sum_{n=0}^{p} \sum_{m=-n}^{n} f_n(r)\widehat{\sigma}_{n,m} Y_n^m, \tag{3.13}$$

with $\widehat{\sigma}_{n,m}$ spherical harmonic coefficients of $\sigma$. Our goal is to evaluate the flux $\boldsymbol{n}^T\nabla\varphi$ at a target point $\boldsymbol{x}$ with spherical coordinates $(r, \theta, \phi)$, where $\boldsymbol{n}(\boldsymbol{x})$ is the target outward unit normal. Let $\varphi_n(r, \theta, \phi) = f_n(r)Y_n^m$. The normal derivative of $\varphi_n$ if given by

$$\boldsymbol{n}^T\nabla\varphi_n = \boldsymbol{n}^T\left(\frac{\partial\varphi_n}{\partial r}\boldsymbol{e}_r + \nabla_\gamma\varphi_n\right), \tag{3.14}$$

$$= f_n'(r)(\boldsymbol{n}^T\boldsymbol{e}_r)Y_n^m + f_n(r)(\boldsymbol{n}^T\nabla_\gamma Y_n^m). \tag{3.15}$$

If the target $\boldsymbol{x}$ is on a sphere centered at the origin, then $\boldsymbol{n}(\boldsymbol{x}) = \boldsymbol{e}_r$, and we need only replace $f_n(r)$ with its derivative in Eq. (3.13) to evaluate $\boldsymbol{n}^T\nabla\varphi$ at $\boldsymbol{x}$. This formula, however, is not valid when $\boldsymbol{n}(\boldsymbol{x})$ has tangential components, as is the case when evaluating the field on one or multiple target spheres. Substitution of Eq. (3.15) yields

$$\boldsymbol{n}^T\nabla\varphi(\boldsymbol{x}) = \sum_{n=0}^{p} \sum_{m=-n}^{n} f_n'(r)(\boldsymbol{n}^T\boldsymbol{e}_r)\widehat{\sigma}_{n,m}Y_n^m + f_n(r)\widehat{\sigma}_{n,m}(\boldsymbol{n}^T\nabla_\gamma Y_n^m) \tag{3.16}$$

We can either evaluate Eq. (3.16) directly, or equivalently use the formula in Eq. (2.9) to write $\boldsymbol{n}^T\nabla_\gamma Y_n^m$ as a linear combination of $Y_{n-1}^m, Y_n^m$ and $Y_{n+1}^m$. Either way, we can observe that the resulting coefficients depend on $(r, \theta, \phi)$, and that the flux does not "diagonalize" in the original spherical harmonic basis.

**Traction calculation (Stokes).** For the velocity field $\boldsymbol{U}$ given by a linear combination of $\mathcal{S}[\boldsymbol{\mu}]$ and $\mathcal{D}[\boldsymbol{\mu}]$, we can for instance evaluate it at any target point on the surface of the unit sphere and its exterior using Eqs. (3.6 - 3.11) to derive a formula of the form

$$\boldsymbol{U}(\boldsymbol{x}) = \sum_{n=0}^{p} \sum_{m=-n}^{n} g_n^V(r)\widehat{\mu}_{n,m}^V\boldsymbol{V}_n^m + (g_n^{VW}(r)\widehat{\mu}_{n,m}^V + g_n^W(r)\widehat{\mu}_{n,m}^W)\boldsymbol{W}_n^m + g_n^X(r)\widehat{\mu}_{n,m}^X\boldsymbol{X}_n^m \tag{3.17}$$

with $\{\widehat{\mu}_{n,m}^V, \widehat{\mu}_{n,m}^W, \widehat{\mu}_{n,m}^X\}$ vector spherical harmonic coefficients of $\boldsymbol{\mu}$. We then want to evaluate the traction of $\boldsymbol{U}$ at a given target point $\boldsymbol{x}$.

We again look at a generic term of Eq. (3.17) $\boldsymbol{u} = g(r)\boldsymbol{Z}$, for $\boldsymbol{Z} = \boldsymbol{V}_n^m, \boldsymbol{W}_n^m$ or $\boldsymbol{X}_n^m$, and let the associated pressure be $p = q(r)Y_n^m$. For the Newton stress tensor $\nabla\boldsymbol{u} + \nabla^T\boldsymbol{u} - pI$, we separate radial and tangential derivatives on the shear stress term $\nabla\boldsymbol{u} + \nabla^T\boldsymbol{u}$:

$$\nabla\boldsymbol{u} + \nabla^T\boldsymbol{u} = \frac{\partial\boldsymbol{u}}{\partial r}\boldsymbol{e}_r^T + \boldsymbol{e}_r\left(\frac{\partial\boldsymbol{u}}{\partial r}\right)^T + \nabla_\gamma\boldsymbol{u} + \nabla_\gamma^T\boldsymbol{u}.$$

The traction force generated by $\boldsymbol{u}$ at a target point is then given by:

$$(\nabla\boldsymbol{u} + \nabla\boldsymbol{u}^T - pI)\boldsymbol{n} = g'(r)(\boldsymbol{e}_r^T\boldsymbol{n}\boldsymbol{Z} + (\boldsymbol{Z}^T\boldsymbol{n})\boldsymbol{e}_r) + g(r)(\nabla_\gamma\boldsymbol{Z} + \nabla_\gamma^T\boldsymbol{Z})\boldsymbol{n} + q(r)Y_n^m\boldsymbol{n}. \tag{3.18}$$

For targets on the sphere, it is again possible to use the fact that $\boldsymbol{n} = \boldsymbol{e}_r$ to simplify Eq. (3.18), and obtain a "diagonal" formula for the traction of $\boldsymbol{u}$ of the form $a_n(r)\boldsymbol{V}_n^m + b_n(r)\boldsymbol{W}_n^m + c_n(r)\boldsymbol{X}_n^m$, as in Eqs 3.6 - 3.8 and 3.9 - 3.11:

$$\mathcal{K}[\boldsymbol{V}_n^m](r,\theta,\phi) = \begin{cases} \dfrac{-2n(n+2)}{(2n+1)(2n+3)}\boldsymbol{V}_n^m r^{-n-3} & r > 1, \\ \dfrac{2n^2+4n+3}{(2n+1)(2n+3)}\boldsymbol{V}_n^m r^n + \dfrac{-(n+1)(n-1)}{2n+1}\boldsymbol{W}_n^m(r^n - r^{n-2}) & r < 1. \end{cases} \tag{3.19}$$

$$\mathcal{K}[\boldsymbol{W}_n^m](r,\theta,\phi) = \begin{cases} \dfrac{-(2n^2+1)}{(2n+1)(2n-1)}\boldsymbol{W}_n^m r^{-n-1} + \dfrac{n(n+2)}{2n+1}\boldsymbol{V}_n^m(r^{-n-1} - r^{-n-3}) & r > 1, \\ \dfrac{2(n+1)(n-1)}{(2n+1)(2n-1)}\boldsymbol{W}_n^m r^{n-2} & r < 1 \end{cases} \tag{3.20}$$

$$\mathcal{K}[\boldsymbol{X}_n^m](r,\theta,\phi) = \begin{cases} \dfrac{-(n+2)}{2n+1}\boldsymbol{X}_n^m r^{-n-2} & r > 1, \\ \dfrac{(n-1)}{2n+1}\boldsymbol{X}_n^m r^{n-1} & r < 1. \end{cases} \tag{3.21}$$

Finally, we must derive a general formula, as $\boldsymbol{n}$ will generally have a nontrivial tangential component. In Fig. 1 we present a spectral analysis in vector spherical harmonics of the functions $(\boldsymbol{e}_r^T\boldsymbol{n})\boldsymbol{Z}+(\boldsymbol{Z}^T\boldsymbol{n})\boldsymbol{e}_r$, $(\nabla_\gamma\boldsymbol{Z} + \nabla_\gamma^T\boldsymbol{Z})\boldsymbol{n}$ and $Y_n^m\boldsymbol{n}$ for $\boldsymbol{Z} = \boldsymbol{V}_n^m, \boldsymbol{W}_n^m, \boldsymbol{X}_n^m$ and $\boldsymbol{n} = \boldsymbol{e}_r, \boldsymbol{e}_\theta, \boldsymbol{e}_\phi$, for $n \leq p = 4$.

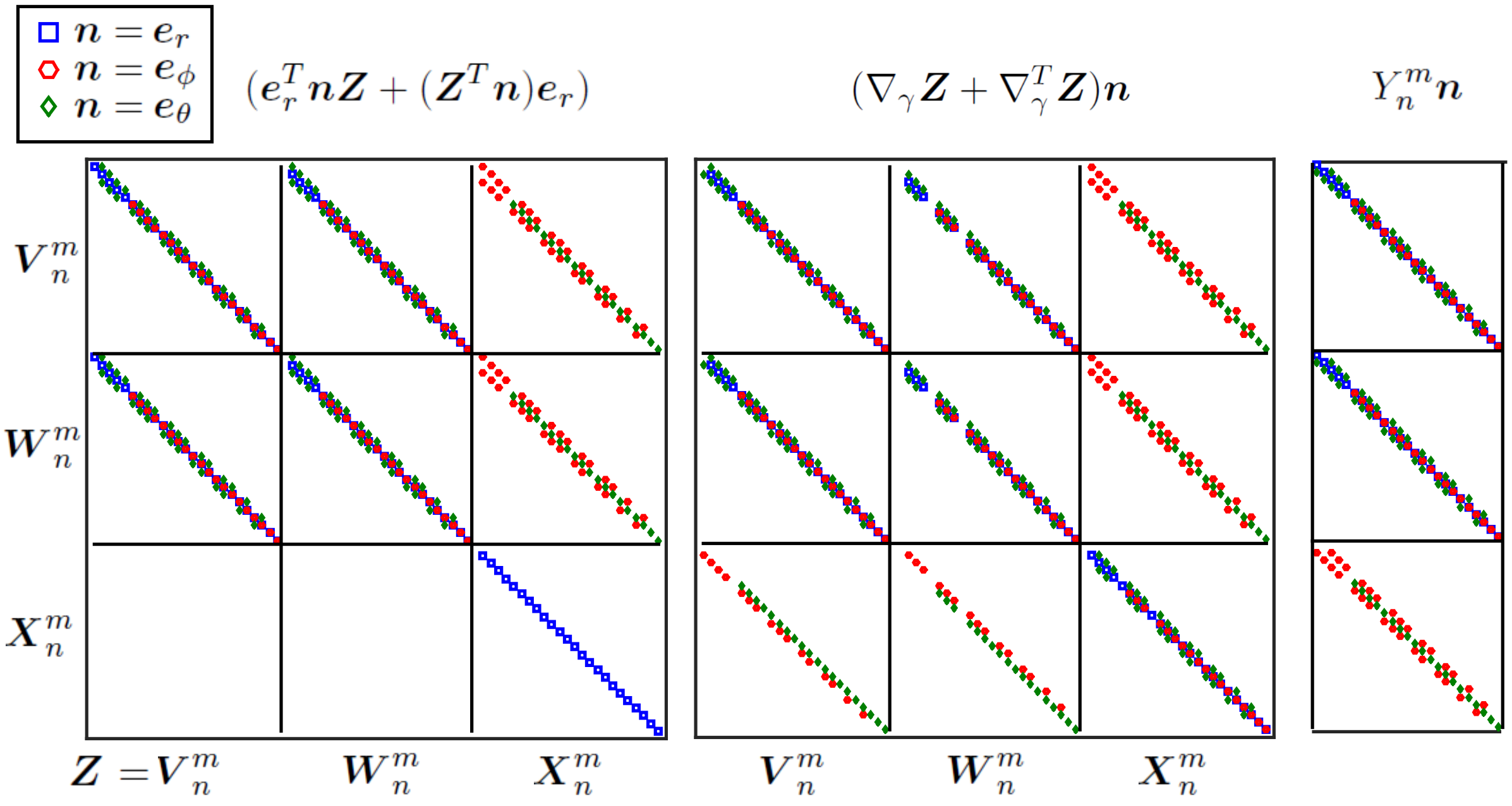

Figure 1: TRACTION COEFFICIENTS. *Coefficients in the spherical harmonic expansion, for* $p = 4$, *of radial, tangential and pressure components of the traction force corresponding to* $\boldsymbol{u} = f(r)Z(\theta, \phi)$ *and outward normal vector* $\boldsymbol{n}$. *Spherical harmonic coefficients for each of the three types are ordered according to* $|m|$ *and then* $n$ *in ascending order. Coefficients for* $\boldsymbol{n} = \boldsymbol{e}_r, \boldsymbol{e}_\theta$ *and* $\boldsymbol{e}_\phi$ *are represented as blue squares, red hexagons and green rhombuses, respectively.*

We observe that the traction generated by the field $\boldsymbol{u} = g(r)\boldsymbol{Z}$ may generally be expanded in terms of spherical harmonics of degree $n-1, n, n+1$ and order $m$. As mentioned previously, the coefficients for $\boldsymbol{n} = \boldsymbol{e}_r$ (in blue on this figure) are only of degree $n$, and so they lie only in the corresponding diagonals.

In order to evaluate the traction of $\boldsymbol{U}$, we take the expansion in Eq. (3.17) and apply the formula Eq. (3.18), with $\boldsymbol{n} = \nu_r \boldsymbol{e}_r + \nu_\theta \boldsymbol{e}_\theta + \nu_\phi \boldsymbol{e}_\phi$. The matrices obtained from analysis in vector spherical harmonics such as in Fig. 1 map the coefficients in each of the resulting 9 terms to an expansion in vector spherical harmonics. In order to accelerate this computation, these matrices are usually precomputed.

# 4 Numerical discretization

While one can compute the solution of the Stokes equations (exterior to spheres) directly using spherical harmonic analysis—as is done in classical multipole methods [8, 37]—the advantage of BIO evaluation and spectral formulas for various kernels derived in the previous section is two-fold. First, we can obtain well-conditioned linear systems capitalizing on the existing work on second-kind BIE formulations for different boundary conditions. Second, we can evaluate the long-range particle-particle hydrodynamic interactions using the standard Stokes FMM [16] instead of developing tailored fast algorithms for different problems.

We now describe a spectrally accurate evaluation scheme for boundary integral operators defined on the spherical particles with boundaries $\{\Gamma_i\}_{i=1}^n$. For a given kernel $K(\boldsymbol{x}, \boldsymbol{y})$, we must evaluate integrals of the form:

$$u(\boldsymbol{x}) = \int_{\Gamma_{\rm src}} K(\boldsymbol{x}, \boldsymbol{y})\sigma(\boldsymbol{y})\, d\Gamma \quad \text{for} \quad \boldsymbol{x} \in \Gamma_{\rm trg}. \tag{4.1}$$

If $\Gamma_{\rm src}$ and $\Gamma_{\rm trg}$ are distinct surfaces and the target $\boldsymbol{x}$ is not close to the integration domain $\Gamma_{\rm src}$, the integrand is a smooth function. As is suggested in our approach, we sample the integral density $\sigma(\boldsymbol{y})$ at points corresponding to a grid in the parametric domain given by

$$\{\theta_j = \cos^{-1}(t_j),\, j = 0, \dots p\} \quad \text{and} \quad \left\{\phi_k = \frac{2\pi k}{2p+2},\, k = 0, \dots, 2p+1\right\}, \tag{4.2}$$

where $t_j$'s are the $(p+1)$-point Gauss-Legendre nodes on $[-1, 1]$. We can then use fast transforms to obtain a truncated expansion in scalar or vector spherical harmonics of order $p$.

## 4.1 Numerical integration schemes

**Smooth integrands.** A standard choice for the numerical integration scheme may thus be obtained using the trapezoidal rule in the azimuthal direction and the Gauss-Legendre quadrature in the polar direction. Then, the following quadrature rule for smooth integrands is spectrally convergent:

$$\int_{\Gamma_{\text{src}}} \sigma(\boldsymbol{y}) d\Gamma = \int_0^{2\pi} \int_0^{\pi} \sigma(\boldsymbol{y}(\theta, \phi)) W(\theta, \phi)\, d\theta d\phi, \tag{4.3}$$

$$= \sum_{j=0}^{p} \sum_{k=0}^{2p+1} \frac{2\pi \lambda_j}{(2p+2)\sin\theta_j} \sigma(\boldsymbol{y}(\theta_j, \phi_k)) W(\theta_j, \phi_k), \tag{4.4}$$

where $\lambda_j$'s are the Gauss-Legendre quadrature weights and $W$ is the corresponding area element for $\Gamma_{\text{src}}$. We will say two boundaries $\Gamma_{\text{src}}$ and $\Gamma_{\text{trg}}$ are well-separated if

$$\text{dist}(\Gamma_{\text{src}}, \Gamma_{\text{trg}}) \geq \eta \max(diam(\Gamma_{\text{src}}), diam(\Gamma_{\text{trg}})), \tag{4.5}$$

where $\eta > 0$ is a heuristic parameter determined based on the user-defined precision. For target points $\boldsymbol{x}$ located in well-separated surfaces, we consider the integrand sufficiently smooth and apply this quadrature rule.

**Singular and near-singular integrands.** When $\Gamma_{\text{trg}} = \Gamma_{\text{src}}$ in (4.1), the integrand becomes singular. On the other hand, when the two surfaces are located close to each other (i.e., (4.5) is not satisfied), it is well-known that the integrand is numerically close to singular and the smooth quadrature rule (4.4) ceases to be effective, as substantial oversampling would be required.

In both cases, we will consider $u(\boldsymbol{x})$ to be the solution to the corresponding exterior problem (e.g. for Stokes or Laplace), for a spherical coordinate system centered at $\boldsymbol{x}^c_{\text{src}}$. Since our kernel $K(\boldsymbol{x}, \boldsymbol{y})$ is in general a combination of single and double layer potentials and their derivatives, we can use the formulas derived in §3 to evaluate $u$ given the spherical harmonic expansion of the source density $\sigma$. For example, in the scalar case, we compute the finite-term approximation of order $p$:

$$\sigma = \sum_{n=0}^{p} \sum_{m=-n}^{n} \widehat{\sigma}_{n,m}\, Y_n^m. \tag{4.6}$$

We store the resulting coefficients in a vector $\widehat{\boldsymbol{\sigma}}$ of length $(p+1)^2$. Evaluating these formulas at target points $\boldsymbol{x}$ with spherical coordinates $\boldsymbol{x} - \boldsymbol{x}^c_{\text{src}} \sim (r, \theta, \phi)$, the solution can be written as:

$$u(\boldsymbol{x}) = \sum_{n=0}^{p} \sum_{m=-n}^{n} \kappa_{n,m}(\boldsymbol{x})\, Y_n^m(\theta, \phi), \tag{4.7}$$

where the vector of coefficients $\boldsymbol{\kappa}(\boldsymbol{x})$ may be computed as $\boldsymbol{\kappa}(\boldsymbol{x}) = \mathcal{F}(\boldsymbol{x})\widehat{\boldsymbol{\sigma}}$, with $\mathcal{F}(\boldsymbol{x})$ a diagonal or banded linear operator.

We note that for in-surface (self-interaction) evaluation, the fact that the operators of interest diagonalize (Theorems 1 and 2) implies that $\mathcal{F}$ does not depend on the target point $\boldsymbol{x}$, and so $\boldsymbol{\kappa}$ may be computed in $O(p^2)$ operations. In this case, $\{\kappa_{n,m}\}$ constitute the spherical harmonic coefficients of order $p$ for $u$ as a function on the sphere. $u(\boldsymbol{x})$ may thus be evaluated using a fast inverse transform, which is $O(p^3 \log p)$, as we choose to evaluate Legendre functions directly, employing the FFT acceleration in the longitudinal direction.

For target points in neighboring surfaces, $\mathcal{F}(\boldsymbol{x})$ is at the very least dependent on $r$, and for examples such as the flux and traction calculations, it is banded. As detailed in §3.3 for calculating the traction, parts of this operator may be precomputed, as those shown in Fig. 1, for efficient evaluation. For a set of $n_{trg}$ target points, computing $\kappa(\boldsymbol{x})$ and evaluating $u(\boldsymbol{x})$ is performed in $O(n_{trg}p^2)$ operations. In the context of evaluation in a suspension of $n$ bodies, this implies a computational cost for direct near-evaluation scaling as $O(np^4)$.

## 4.2 Fast near-singular evaluation

For all BIOs of interest, particularly for the case of the Stokes traction kernel, direct evaluation of the near-singular integration formulas in §3 can become quite expensive for moderate to high order $p$. We present an alternate fast algorithm in the context of their evaluation on a set of polydisperse neighbor target spheres, each sampled at the corresponding grid as in Eq. (4.2). For ease of presentation, we focus on evaluation of exterior problems with no bounding geometry (e.g. all target spheres lie outside the source sphere). Note, however, that the methods presented below extend in a straightforward manner to the case in which the target lies inside the source sphere.

The algorithm below relies on the fact that key computations simplify when the north poles for source and target spheres are aligned, and so we may couple it with fast routines for rotation of spherical harmonics expansions. This allows us to perform near-singular evaluation identical to that of the $O(p^3 \log p)$ in-surface evaluation. We present an experimental comparison for increasing order $p$, demonstrating considerable speed-ups against direct evaluation for both scalar and vector expansions.

This method may be interpreted as a variant of well-established "point-and-shoot" strategies used to accelerate the translation of multipole expansions for wideband Fast Multipole Methods for the Helmholtz equation in three dimensions [6, 7]. We end this section with a brief discussion of the scalar case as well as its extension to the Stokes equations.

### FFT-accelerated near-evaluation

As indicated in §4.1, the two reasons we are able to perform fast in-surface evaluation of Eq. (4.7) are that coefficients $\kappa_{n,m}$ are independent of the target point $\boldsymbol{x}$, and that the fast inverse transform then arranges computation in two stages: direct evaluation of Legendre associated functions $\{P_n^{|m|}(cos\theta_j)\}_{j=0}^{p}$ and fast evaluation along each latitude $\theta_j$ using an FFT. In other words, we rearrange the sums in Eq. (4.7):

$$u(\boldsymbol{x}) = \sum_{m=-p}^{p} \left( \sum_{n=|m|}^{p} \widetilde{\kappa}_{n,m} P_n^{|m|}(\cos\theta) \right) e^{im\phi} = \sum_{m=-p}^{p} G_m(\cos\theta) e^{im\phi}, \tag{4.8}$$

with $\widetilde{\kappa}_{n,m}$ equal to the coefficients times the constant factor in $Y_n^m$.

Replicating this strategy for targets in a neighboring sphere generally fails on both counts: $\kappa_{n,m}$ are target-dependent and the angles $(\theta, \phi)$ at the target sphere are generally not equispaced along (or aligned with) lines of constant latitude or longitude.

Consider, however, the special case in which source and target sphere north poles are aligned; that is, the target's center $C$ is of the form $[0\ 0\ C_z]$ for $C_z > 0$ and we assume its computational grid is a translation of the source grid. As can be seen in Fig. 2(2), target points $\boldsymbol{x}_{j,k}$ can be grouped in $p+1$ discs parallel to the $xy$-plane, with spherical coordinates $\boldsymbol{x}_{j,k} = (r_j, \theta_j, \phi_k)$ where $\phi_k$ are equispaced and $(r_j, \theta_j)_{j=0}^{p}$ are constant on each disc. The FFT-accelerated algorithm detailed above is thus viable as long as $\kappa_{n,m}$ are *only* dependent on $(r, \theta)$.

This can be readily observed to be the case for $\mathcal{S}^L$ and $\mathcal{D}^L$ in §3.1, as the coefficients in Eq. (3.3) and Eq. (3.4) depend only on $r$. For the flux calculation in §3.3, coefficients in Eq. (3.16) ultimately are a linear combination of terms of the form $\{f_n'(r)\nu_r, f_n(r)\nu_\theta, f_n(r)\nu_\phi\}$, with $(\nu_r, \nu_\theta, \nu_\phi)$ spherical coordinate coefficients of the normal vector $\boldsymbol{n}(\boldsymbol{x})$.

Fortunately, in the pole-aligned setting, these coordinates are independent of $\phi$. For a target sphere of radius $R$, we have

$$\begin{bmatrix} \nu_r(\theta) & \nu_\theta(\theta) & \nu_\phi(\theta) \end{bmatrix} = \begin{bmatrix} \frac{R + C_z \cos\theta}{\sqrt{R^2 + 2C_z \cos\theta + C_z^2}} & \frac{C_z \sin\theta}{\sqrt{R^2 + 2C_z \cos\theta + C_z^2}} & 0 \end{bmatrix}. \tag{4.9}$$

This in turn simplifies Eq. (3.16):

$$\boldsymbol{n}^T \nabla \varphi(\boldsymbol{x}) = \sum_{n=0}^{p} \sum_{m=-n}^{n} f_n'(r)\nu_r(\theta)\widehat{\sigma}_{n,m} Y_n^m + f_n(r)\nu_\theta(\theta)\widehat{\sigma}_{n,m} \frac{\partial Y_n^m}{\partial \theta} \tag{4.10}$$

and thus the conditions for applying our FFT-accelerated method are met (computing $G_m(\cos\theta)$ in fact simply involves a linear combination of $P_n^{|m|}$ and its derivative).

**General case:** The previous result suggests a general algorithm in three stages as depicted in Fig. 2. We assume a source sphere of unit radius centered at the origin, and a target sphere with center $C$.

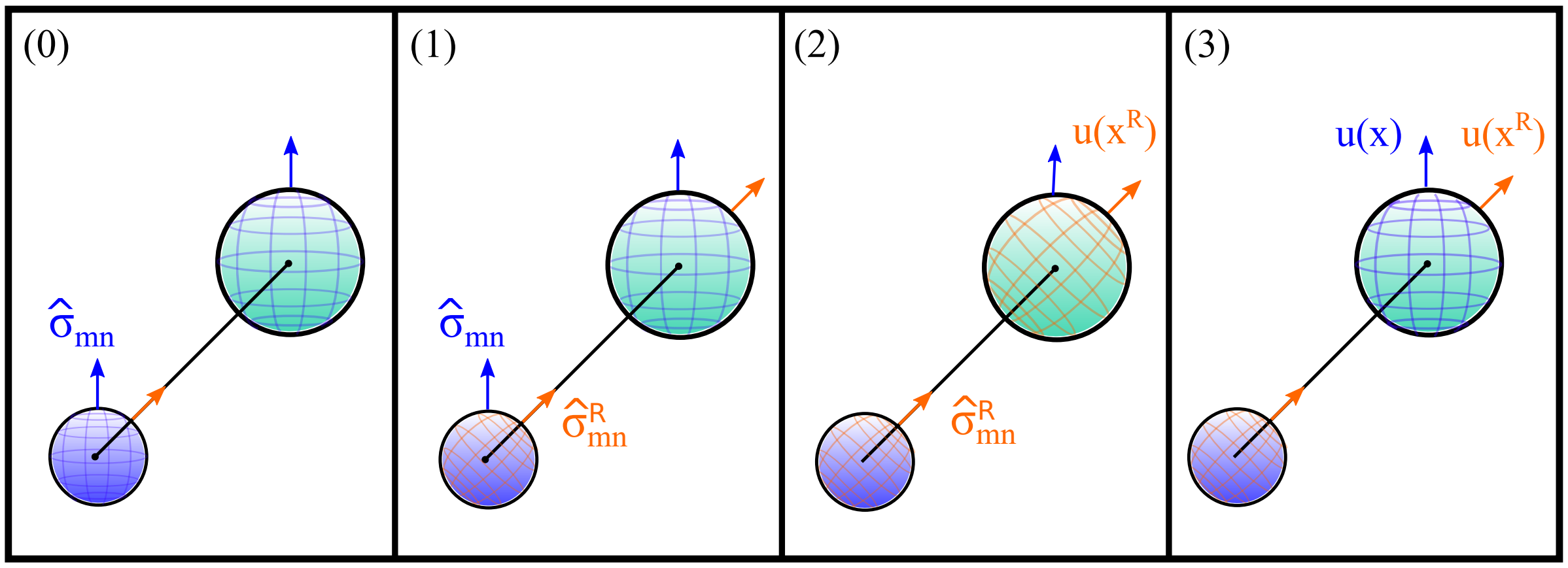

Figure 2: FFT-accelerated near evaluation algorithm.
(0) Input: *spherical harmonic coefficients* $\widehat{\sigma}_{n,m}$ *for density* $\sigma$ *defined on the original grid (blue),*
(1) Density rotation: *We rotate the computational grid so that the north pole points in the direction of target center* $C$*. We compute equivalent coefficients* $\widehat{\sigma}^R_{n,m}$ *for the rotated grid (orange),*
(2) Fast evaluation at translated grid: *We apply our fast near-evaluation algorithm, computing function evaluations* $u(\boldsymbol{x}^R)$ *at a target grid aligned with the rotated pole (both orange),*
(3) Rotate function samples: *We rotate the target grid to match the original target grid pole (blue). Given samples* $u(\boldsymbol{x}^R)$*, we compute* $u(\boldsymbol{x})$ *at original target points* $\boldsymbol{x}$ *(blue).*

Rotation of spherical harmonic expansions to a new spherical grid can be accomplished in $O(p^3)$ operations using Wigner rotation matrices [48], see [17] for a detailed discussion. Stage 1 consists of one such rotation. In stage 3, we use a fast forward transform to compute coefficients $\widehat{u}^R_{n,m}$, rotate to the desired target pole, and evaluate at target points with one fast inverse transform. This algorithm thus has an overall cost of $O(p^3 \log p)$.

We note that when evaluating self and near-field interactions, it is advantageous to delay the final evaluation of $u(x)$ in stage (3), storing $\widehat{u}_{n,m}$ instead. For each target sphere, we then add spherical harmonic coefficients from itself and all its neighbors, and then perform one fast inverse transform to evaluate.

### Vector spherical harmonics FFT-accelerated algorithm

The fast vector spherical harmonic transforms employed in this work rely on three scalar transforms via sparse (pentadiagonal) transformations to and from expansions of a field in spherical coordinates. We can then generalize the FFT-accelerated algorithm presented above to evaluate expressions of the form Eq. (3.17) as long as the resulting coefficients in each coordinate direction require $O(p^3)$ to compute, and depend only on $(r, \theta)$ at target points when its pole is aligned with that of the computational grid.

This may again be readily concluded for the Stokes single and double layer potentials: all coefficients depend only on $r$, and we need only apply the transformations in Eqs. (2.16 - 2.17) and 2.15 (one for each of the $(p+1)$ discs, for a total of $O(p^3)$ operations), followed by one application of the FFT-accelerated evaluation algorithm for each spherical coordinate field.

For the traction kernel calculation in §3.3, as it was the case for the flux for Laplace, coefficients for the traction of $\boldsymbol{U}$ in the basis $\{\boldsymbol{V}_n^m, \boldsymbol{W}_n^m, \boldsymbol{X}_n^m\}$ are ultimately linear combinations of functions $g_n(r)$ multiplied by $\nu_r, \nu_\theta$ and $\nu_\phi$, obtained by the application of the 3 banded matrices depicted in Fig. 1. Since $\nu_r$ and $\nu_\theta$ depend only on $\theta$ and $\nu_\phi = 0$, the resulting coefficients in the vector spherical harmonic bases depend only on $(r, \theta)$, and we may again apply the FFT-accelerated algorithm. Furthermore, we note that once again, the fact that $\boldsymbol{n}$ is orthogonal to $e_\phi$ removes some of the computation involved, as the traction coefficients in that direction are also zero.

For the general algorithm, we note rotation of spherical harmonic densities and function values can again be readily be computed using scalar transforms.

### Comparison with direct near-evaluation

We conduct an experiment in order to compare serial implementations of the FFT-accelerated algorithm against direct evaluation. We compute expansions for the Laplace and Stokes single layer potentials on 10 random target spheres for increasing order $p$, and measure average timings for both methods and assess their experimental scaling.

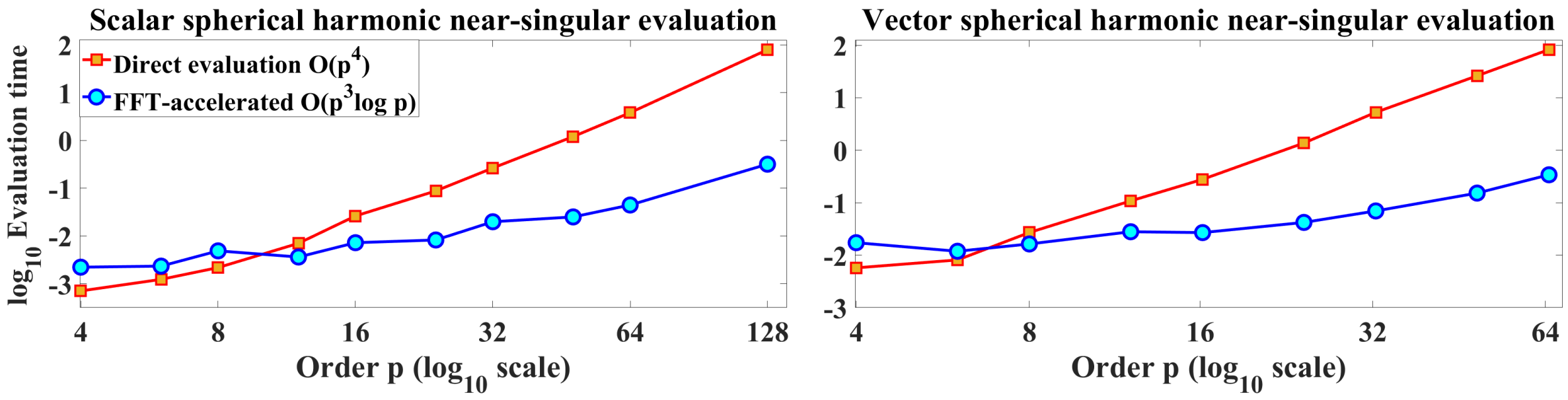

Figure 3: Near-singular evaluation algorithm comparison. *Left: log-log plot comparison of average timings for near-singular evaluation algorithms applied to the Laplace single layer potential on a random target sphere. Right: analogous comparison for vector spherical harmonics algorithms applied to the valuation of the Stokes single layer potential.*

In Fig. 3, we can observe that although both methods are comparably fast for $p \leq 8$, the FFT-accelerated approach provides a significant speed-up for moderate and high $p$, achieving a two orders of magnitude acceleration for $p \geq 64$. We note that experimental scaling shows the FFT-based method is considerably efficient, as it only gets close to $O(p^3)$ for $p \geq 48$.

### Translation operators for solid spherical harmonics

As indicated in §3.1, the solutions $f_n(r)Y_n^m$ to the Laplace equation we use in the expansions of our layer potentials are known as *solid spherical harmonics*; distinguishing between solutions valid in the interior (*regular solid harmonics*, denoted $R_n^m(\boldsymbol{x})$) and those valid in the exterior (*irregular solid harmonics*, denoted $I_n^m(\boldsymbol{x})$) with $\boldsymbol{x} = (r, \theta, \phi)$.

The "point-and-shoot" algorithms in wideband fast multipole methods employ a strategy based on the addition theorems for solid harmonics [40] in lieu of the FFT acceleration in stage 2 (Fig. 2), in order to translate wave expansions from children to parent (M2M) and viceversa (L2L). For a detailed derivation we refer the reader to Chapter 5 of [7] and [19]. The key result (Thm. 3 in [6]) is that when poles for target and source spheres are aligned, translation operators simplify considerably, yielding expressions of the form:

$$u(\boldsymbol{x}) = \sum_{n=0}^{p} \sum_{m=-n}^{n} \kappa_{n,m} f_n(r) Y_n^m(\theta, \phi) \tag{4.11}$$

$$= \sum_{\lambda=0}^{p} \sum_{m=-\lambda}^{\lambda} \left( \sum_{n=0}^{p} \kappa_{n,m} B_{\lambda,n,m} \frac{r_T^{\lambda}}{C_z^{\lambda+n+1}} \right) Y_\lambda^m(\theta_T, \phi_T) \tag{4.12}$$

where $(r_T, \theta_T, \phi_T)$ are coordinates for $\boldsymbol{x}$ in the target sphere. This allows us to translate the full expansion to the target grid with spherical coordinates in $O(p^3)$ operations. We can apply these formulas directly to single and double layer potentials, and translate the relevant differential operators to compute any derivatives of interest on the target grid. Derivatives in the radial direction can be computed analytically, and spectral methods can be used for $(\theta_T, \phi_T)$. Finally, we note that these methods can be extended to vector spherical harmonics for the Stokes equations employing the appropriate solid harmonic addition theorems [12, 13].

### 4.3 Fast integral evaluation algorithm

In order to produce a fast, optimal complexity evaluation scheme for integral operators in this setting, we use the numerical integration schemes described above within the framework of the FMM. For both scalar and vector case, the total number of degrees of freedom in these systems is $N = O(np^2)$.

We use the point-FMM 3D libraries for Laplace and Stokes potentials developed by Gimbutas & Greengard [15, 16], evaluating far-interactions with the smooth quadrature in Eq. (4.4). Self-interactions are computed using our singular-evaluation scheme. For neighboring surfaces which are not sufficiently separated according to Eq. (4.5), we make the necessary correction using the near-singular scheme. Overall, this yields a spectrally-accurate evaluation scheme with complexity $O(np^3 \log p)$. In §6 we confirm this scaling experimentally, and we apply it to evaluate various integral operators relevant to the study of suspensions of rigid bodies in Stokesian flows.

We note that since state-of-the-art three dimensional FMM codes rely on spherical harmonic representation for multipole and local series expansions, our evaluation scheme may also be used in the context of an FMM algorithm whose input and output are spherical harmonic coefficients of degree $p$ on each surface $\Gamma_i$. Since this is an optimization that only affects the source to multipole operator, we do not expect it to drastically change overall performance.

## 5 BIE formulations and case studies

In this section, we discuss different classes of problems encountered in applications, a few are depicted in Figure 4, and demonstrate the new solution procedure. We then perform a series of experiments to verify the convergence and scaling properties of our algorithm in Section 6.

In all of the case studies presented next, we assume that $n_b$ rigid spherical particles $\{D_i\}_{i=1}^{n_b}$ are suspended in a Stokes flow with viscosity $\nu = 1$ and in free-space. For a given velocity field $\boldsymbol{u}(\boldsymbol{x}) \in \mathbb{R}^3$ at an arbitrary point $\boldsymbol{x}$ in the fluid domain, we denote the corresponding fluid pressure and stress tensors by $p$ and $\sigma$, respectively. The governing equations are then given by:

$$-\Delta \boldsymbol{u} + \boldsymbol{\nabla} p = \boldsymbol{0}, \quad \boldsymbol{\nabla} \cdot \boldsymbol{u} = 0 \quad \forall \quad \boldsymbol{x} \in \mathbb{R}^3 \setminus \cup_i D_i \tag{5.1}$$

$$\boldsymbol{u}(\boldsymbol{x}) \to \boldsymbol{0} \ \text{ as } \ |\boldsymbol{x}| \to \infty \tag{5.2}$$

Each of the problems discussed next solve the same equations as above but with different boundary and kinematic conditions. While a variety of boundary integral formulations are available, we will present our methods of choice only; [35] discusses several other classical formulations.

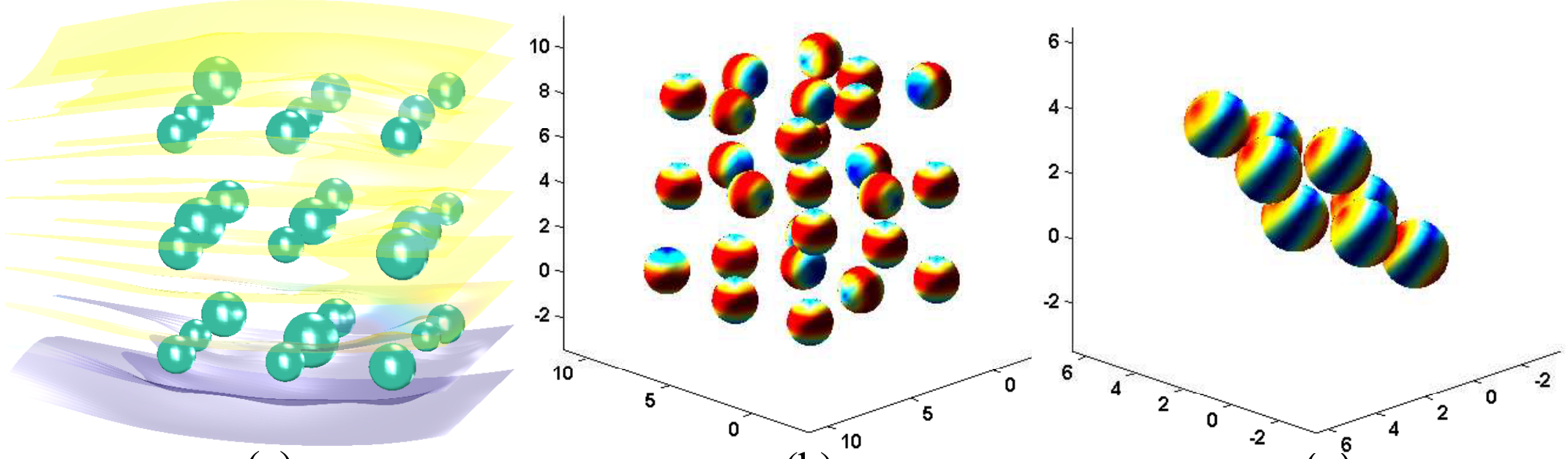

(a) (b) (c)

Figure 4: (a) Porous media flow: constant flow (right to left) past a regular, polydisperse lattice. (b) Squirmer flow: randomly oriented ''pusher" squirmers moving away from a cubic lattice configuration. Sphere surface color is proportional to slip velocity. (c) MHD flow: final configuration for paramagnetic beads after subjecting them to a constant magnetic field $\boldsymbol{H}_0 = [10\ 10\ 10]$. Sphere surface color is proportional to the applied magnetic traction.

**Problem 1: Porous media flow**

In this problem, the particles are static and an exterior flow field $\boldsymbol{u}_\infty(\boldsymbol{x})$ is imposed. Assume that the flow field is given by,

$$\textbf{ansatz:} \qquad \boldsymbol{u}(\boldsymbol{x}) = \boldsymbol{u}_\infty(\boldsymbol{x}) + \Sigma_{k=1}^{n_b}(\mathcal{S}_k + \mathcal{D}_k)[\boldsymbol{\mu}](\boldsymbol{x}) \tag{5.3}$$

Then, taking the limit as $\boldsymbol{x}$ approaches a particle boundary from the exterior and applying the no-slip boundary condition, we get the following

$$\textbf{BIE:} \qquad \left(\frac{1}{2}\boldsymbol{I} + \Sigma_{k=1}^{n_b}(\mathcal{S}_k + \mathcal{D}_k)\right)[\boldsymbol{\mu}](\boldsymbol{x}) = -\boldsymbol{u}_\infty(\boldsymbol{x}) \tag{5.4}$$

Solving the above equation for the unknown density function $\boldsymbol{\mu}$, we can substitute it in (5.3) and obtain the velocity field at any point in the fluid domain.

**Problem 2: Mobility problem**

Here, we are given forces and torques applied on the rigid particles and we need to find the resulting rigid body velocity fields. Let $(\boldsymbol{F}_k, \boldsymbol{T}_k)$ be the pairs of forces and torques applied to the particle with boundary $\Gamma_k$ and let $(\boldsymbol{v}_k, \boldsymbol{\omega}_k)$ be the particle's rigid-body translational and rotational velocities. Then, we have the following conditions on the particle boundary:

$$\boldsymbol{u}(\boldsymbol{x}) = \boldsymbol{v}_k + \boldsymbol{\omega}_k \times (\boldsymbol{x} - \boldsymbol{x}_k^c) \quad \forall\ \boldsymbol{x} \in \Gamma_k, \tag{5.5}$$

$$\int_{\Gamma_k} \boldsymbol{f}\, d\Gamma = \int_{\Gamma_k} \sigma \cdot \boldsymbol{n}\, d\Gamma = -\boldsymbol{F}_k, \quad \int_{\Gamma_k} (\boldsymbol{x} - \boldsymbol{x}_k^c) \times \boldsymbol{f}\, d\Gamma = -\boldsymbol{T}_k, \tag{5.6}$$

where $\boldsymbol{x}^c$ is the centroid and $\boldsymbol{f}$ is the traction. We follow [9, 36] and formulate the integral equations as follows.

$$\textbf{Ansatz:} \qquad \boldsymbol{u}(\boldsymbol{x}) = \Sigma_{k=1}^{n_b}\mathcal{S}_k[\boldsymbol{\mu} + \boldsymbol{\rho}](\boldsymbol{x}). \tag{5.7}$$

The density function $\boldsymbol{\rho}$ is computed directly from the given forces and torques (see [9]). Then, prescribing the total internal stress from the inside of the each particle to be zero (which corresponds to a rigid

body motion), we get

$$\textbf{BIE:} \qquad \left(\frac{1}{2}\boldsymbol{I} + \Sigma_{k=1}^{n_b}\left(\mathcal{K}_k + \mathcal{L}_k\right)\right)[\boldsymbol{\mu}](\boldsymbol{x}) = -\left(\frac{1}{2}\boldsymbol{I} + \Sigma_{k=1}^{n}\mathcal{K}_k\right)[\boldsymbol{\rho}](\boldsymbol{x}), \tag{5.8}$$

where the operator $\mathcal{L}_k[\boldsymbol{\mu}](\boldsymbol{x}) = \int_{\Gamma_k}\boldsymbol{\mu}(\boldsymbol{y})d\Gamma + \left(\int_{\Gamma_k}(\boldsymbol{y}-\boldsymbol{x}_i^c)\times\boldsymbol{\mu}(\boldsymbol{y})d\Gamma\right)\times(\boldsymbol{x}-\boldsymbol{x}_i^c)$ is added to the formulation in order to eliminate the $6n$ dimensional nullspace of the operator $\frac{1}{2}I+\mathcal{K}$ (corresponding to rotations and translations) and select the unique $\boldsymbol{\mu}$ that does not add net forces and torques.

### Problem 3: Multiphysics—MHD flow

We consider an integral equation formulation for the evolution of suspensions of paramagnetic beads in Stokesian fluid, under the application of a constant magnetic field, as presented in [9]. Assuming the fluid is insusceptible and an absence of free currents simplifies the Maxwell equations, allowing for a representation of the magnetic field as the negative gradient of a scalar potential $\phi$. This is a well-known simplification for the magnetostatic case, reducing the Maxwell equations to a Laplace equation for $\phi$ with prescribed jump conditions at the surface boundaries:

$$\Delta\phi(\boldsymbol{x}) = 0 \quad \forall\ \boldsymbol{x} \notin \Gamma, \tag{5.9}$$

$$[[\phi]]_\Gamma = 0, \qquad \left[\left[\mu\frac{\partial\phi}{\partial\boldsymbol{r}}\right]\right]_\Gamma = 0, \tag{5.10}$$

$$\phi \to -\boldsymbol{H_0}\cdot\boldsymbol{r} \ \text{ as } \ |\boldsymbol{x}|\to\infty. \tag{5.11}$$

Here, $\mu$ is the magnetic permeability and $\boldsymbol{H}_0$ is the imposed magnetic field. We formulate the integral equation as follows.

$$\textbf{Ansatz:} \qquad \phi(\boldsymbol{x}) = -\boldsymbol{H}_0\cdot\boldsymbol{x} + \mathcal{S}^L[q](\boldsymbol{x}). \tag{5.12}$$

Applying the boundary conditions and using the standard jump conditions [25, 35] gives us,

$$\textbf{BIE:} \qquad \left(\frac{1}{2}I + \eta\mathcal{K}_\Gamma^L\right)[q](\boldsymbol{x}) = \eta\boldsymbol{H}_0\cdot\boldsymbol{n}, \tag{5.13}$$

where $\mathcal{K}^L$ is the normal derivative of the single layer potential and $\eta = \frac{\mu-\mu_0}{\mu+\mu_0}$.

Being that the fluid is insusceptible, the only coupling between the magnetic field and the flow occurs through the traction forces and torques applied to each particle surface. Therefore the above formulation must be combined with one for the Stokes mobility problem. For a given configuration of bodies with boundaries $\Gamma_k$, we can obtain the corresponding scalar potential $\phi$, and set the incident force field density $\boldsymbol{\rho}$ in Eq. (5.8) to the corresponding magnetic traction on surface $\Gamma_k$ (computed using the Maxwell stress tensor).

### Problem 4: Active particles

A classical model for an active particle in Stokes flow is the so called *squirmer* model [3, 26, 28]. Each particle has a prescribed slip velocity $\boldsymbol{u}_s$, leading to the conditions:

$$\boldsymbol{u}(\boldsymbol{x}) = \boldsymbol{u}_s(\boldsymbol{x}) + \boldsymbol{v}_k + \boldsymbol{\omega}_k\times(\boldsymbol{x}-\boldsymbol{x}_k^c) \quad \forall \quad \boldsymbol{x}\in\Gamma_k, \tag{5.14}$$

$$\int_{\Gamma_k}\boldsymbol{f}\,d\Gamma = \boldsymbol{0} \quad \text{and} \quad \int_{\Gamma_k}(\boldsymbol{x}-\boldsymbol{x}_k^c)\times\boldsymbol{f}\,d\Gamma = \boldsymbol{0}. \tag{5.15}$$

We use the following integral equation formulation.

$$\textbf{Ansatz:} \qquad \boldsymbol{u}(\boldsymbol{x}) = \Sigma_{k=1}^{n_b}\left(\mathcal{S}_k + \mathcal{D}_k\right)[\boldsymbol{\mu}](\boldsymbol{x}). \tag{5.16}$$

Again, taking the limit as $\boldsymbol{x}$ approaches a particle boundary $\Gamma_i$ from the exterior, we get

$$\textbf{BIE:} \qquad \left(\frac{1}{2}\boldsymbol{I} + \Sigma_{k=1}^{n_b}\left(\mathcal{S}_k + \mathcal{D}_k\right)\right)[\boldsymbol{\mu}](\boldsymbol{x}) = \boldsymbol{u}_s(\boldsymbol{x}) + \boldsymbol{v}_i + \boldsymbol{\omega}_i \times (\boldsymbol{x} - \boldsymbol{x}_i^c) \quad \text{for} \quad \boldsymbol{x} \in \Gamma_i. \tag{5.17}$$

The unknowns at every time-step $\{\boldsymbol{v}, \boldsymbol{\omega}, \boldsymbol{\mu}\}$ are computed by solving the above equation coupled with (5.15). Note that the double layer potential doesn't contribute to forces and torques, which simplifies application of (5.15). The particle configurations are then updated using the thus obtained rigid body velocity fields.

**Problem 5: Resistance problem**

The resistance problem is in a sense the inverse of the mobility problem: we are given the particle's rigid body translational and rotational velocities $(\boldsymbol{v}_k, \boldsymbol{\omega}_k)$, and we need to find the corresponding forces and torques $(\boldsymbol{F}_k, \boldsymbol{T}_k)$. Conditions at the boundary Eq. (5.5) and Eq. (5.6) are again imposed.

We apply the classical *completed double-layer* formulation (Power and Miranda [34]) for this problem, which proceeds as follows.

$$\textbf{Ansatz:} \qquad \boldsymbol{u} = \Sigma_{k=1}^{n_b}\left(\mathcal{D}_k + \mathcal{N}_k\right)[\boldsymbol{\psi}](\boldsymbol{x}). \tag{5.18}$$

$$\textbf{BIE:} \qquad \left(\frac{1}{2}\boldsymbol{I} + \Sigma_{k=1}^{n_b}\left(\mathcal{D}_k + \mathcal{N}_k\right)\right)[\boldsymbol{\psi}](\boldsymbol{x}) = \boldsymbol{v}_i + \boldsymbol{\omega}_i \times (\boldsymbol{x} - \boldsymbol{x}_i^c), \tag{5.19}$$

where $\mathcal{N}_k$ is sometimes referred to as a completion flow. A standard choice for this formulation is to add the flow generated by a point force (Stokeslet) and point torque (Rotlet) located at the center of each body $D_k$.

For all the problems presented above, notice that the spectra and the evaluation formulas given in Section 3 are sufficient for their solution.

# 6 Numerical Results

We now conduct a series of tests to validate the accuracy, convergence properties and scaling of numerical integration schemes for the boundary integral operators as described in §4.

**Singular and near-singular integration**

First, we focus on testing the behavior of the near-singular integration scheme described in §4, as compared to that of the smooth quadrature Eq. (4.4) when the target points are brought closer to the source boundary. For this purpose, we first randomly generate a density $\boldsymbol{\sigma}$ on the boundary of the unit sphere:

$$\boldsymbol{\sigma} = \sum_{n=0}^{p_{max}} \sum_{m=-n}^{n} f_n^m \boldsymbol{V}_n^m + g_n^m \boldsymbol{W}_n^m + h_n^m \boldsymbol{X}_n^m$$

for $p_{max} = 16$, imposing an algebraic decay of the coefficients $f, g, h$ with respect to $n$, as is characteristic of smooth data. For each integral operator, we then test performance of both schemes for spherical harmonic expansions of order $p \in 4, 8, 16$. We evaluate the layer potentials on a set of spherical shells of target points at distances $\{1, 10^{-0.5}, \ldots, 10^{-6}\}$ from the surface, and compare with the exact solutions in §3.2.

From the box plots presented in Fig. 5, we can readily observe that the near-singular evaluation scheme is spectrally convergent, and that it remains accurate for target points arbitrarily close to the sphere surface. The smooth quadrature, as expected, is spectrally convergent in the far-field yet it becomes ineffective as our target points approach the surface.

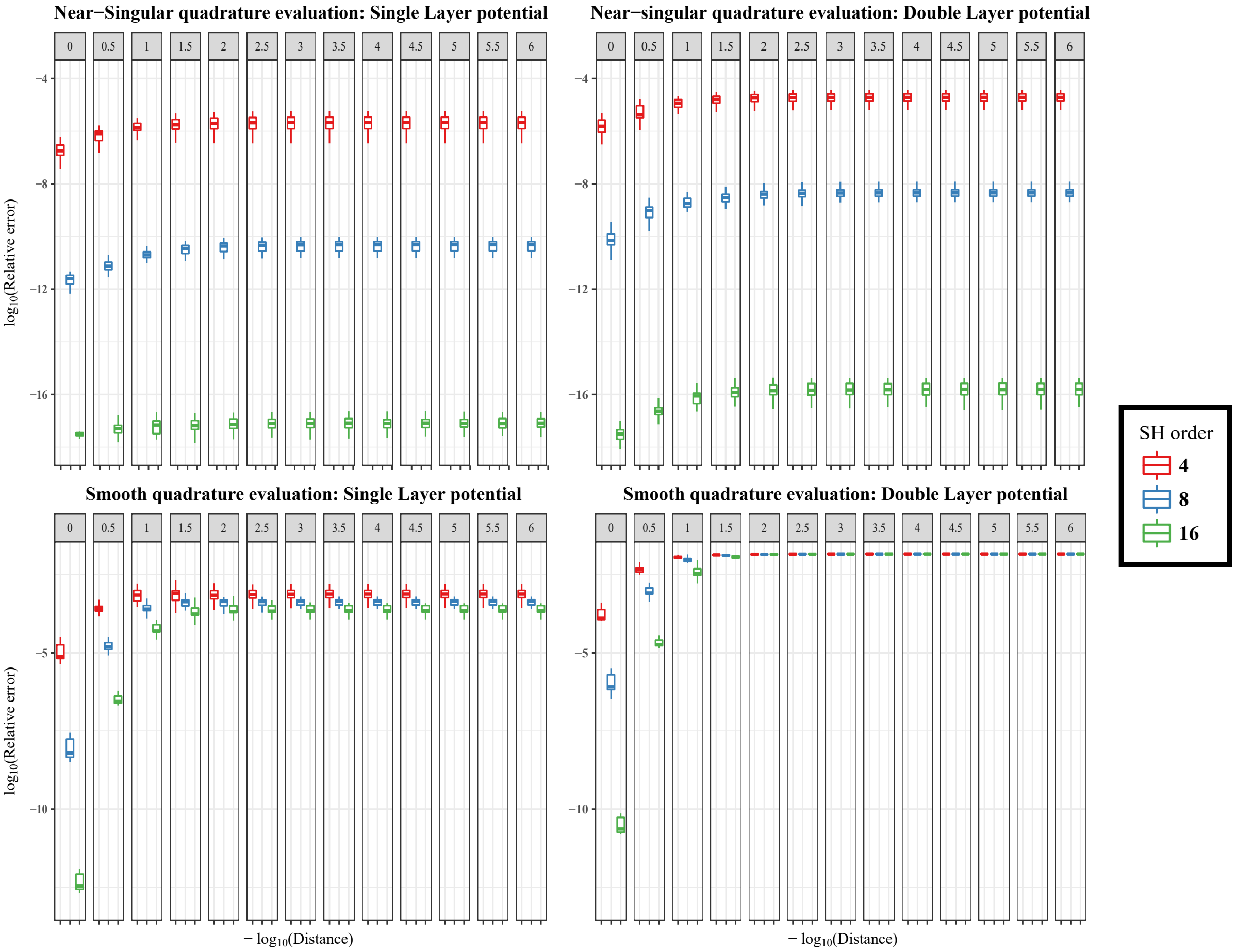

Figure 5: Quadrature performance on synthetic data. *We measure the* log *relative error for the evaluation of* $\mathcal{S}$ *and* $\mathcal{D}$ *for the Stokes equation for both integration schemes. For spherical harmonic order* $p = 4, 8, 16$*, we present box plots categorized by* $-\log$ *of the distance to the unit sphere surface.*

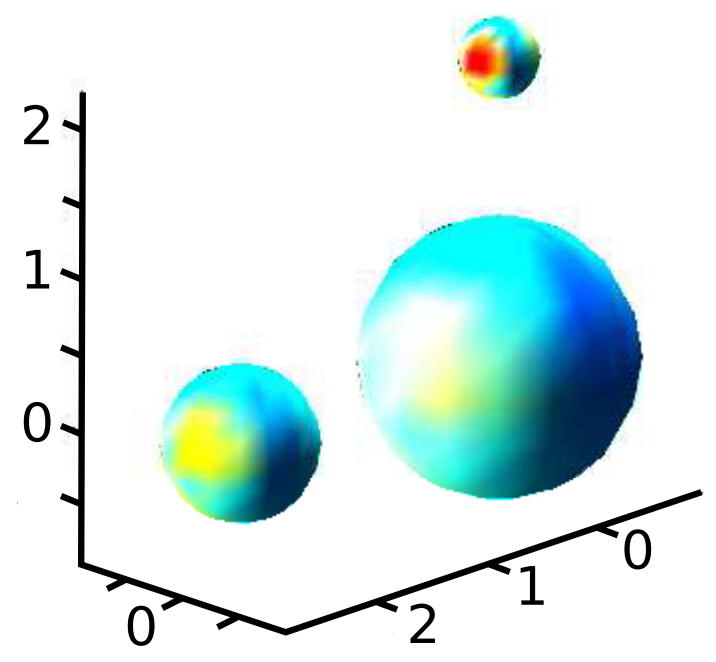

Figure 6: Sphere configuration for BIE data example.

In order to further validate our evaluation scheme, we also test its performance in the solution of boundary integral equations such as those presented on the case studies in §5. We begin with an arrangement of $n_b = 3$ spheres of varying radii (polydisperse) $r_\ell = (0.903, 0.510, 0.262)$ (Fig. 6).

We then pick a random point $\boldsymbol{p}_\ell$ inside each sphere ($||\boldsymbol{p}_\ell - \boldsymbol{x}^c_\ell|| = 0.5 r_\ell$), and produce the flows

generated by point forces (Stokeslets) $\gamma_\ell(x) = \boldsymbol{G}(\boldsymbol{x}, \boldsymbol{p}_\ell)$ and corresponding tractions (Stresslet) $\tau_\ell(x) = \boldsymbol{T}_{i,j,k}(\boldsymbol{x}, \boldsymbol{p}_\ell)\boldsymbol{n}_k(\boldsymbol{x})$. We solve the following boundary integral equations

$$\left(\frac{1}{2}I + \mathcal{K}_\Gamma + \mathcal{L}\right)[\boldsymbol{\mu}](x) = \sum_{\ell=1}^{n_b} \tau_\ell(x) \tag{6.1}$$

$$\left(\frac{1}{2}I + \mathcal{D}_\Gamma + \mathcal{N}\right)[\boldsymbol{\psi}](x) = \sum_{\ell=1}^{n_b} \gamma_\ell(x) \tag{6.2}$$

taken from the mobility (Eq. (5.8)) and resistance (Eq. (5.19)) formulations, respectively. We again set a series of target points around each sphere, at distances $r_\ell\{1, 10^{-0.5}, \ldots, 10^{-6}\}$ from their surface, and measure evaluation error for $\mathcal{S}[\boldsymbol{\mu}]$,$\mathcal{K}[\boldsymbol{\mu}]$ and $\mathcal{D}[\boldsymbol{\psi}]$.

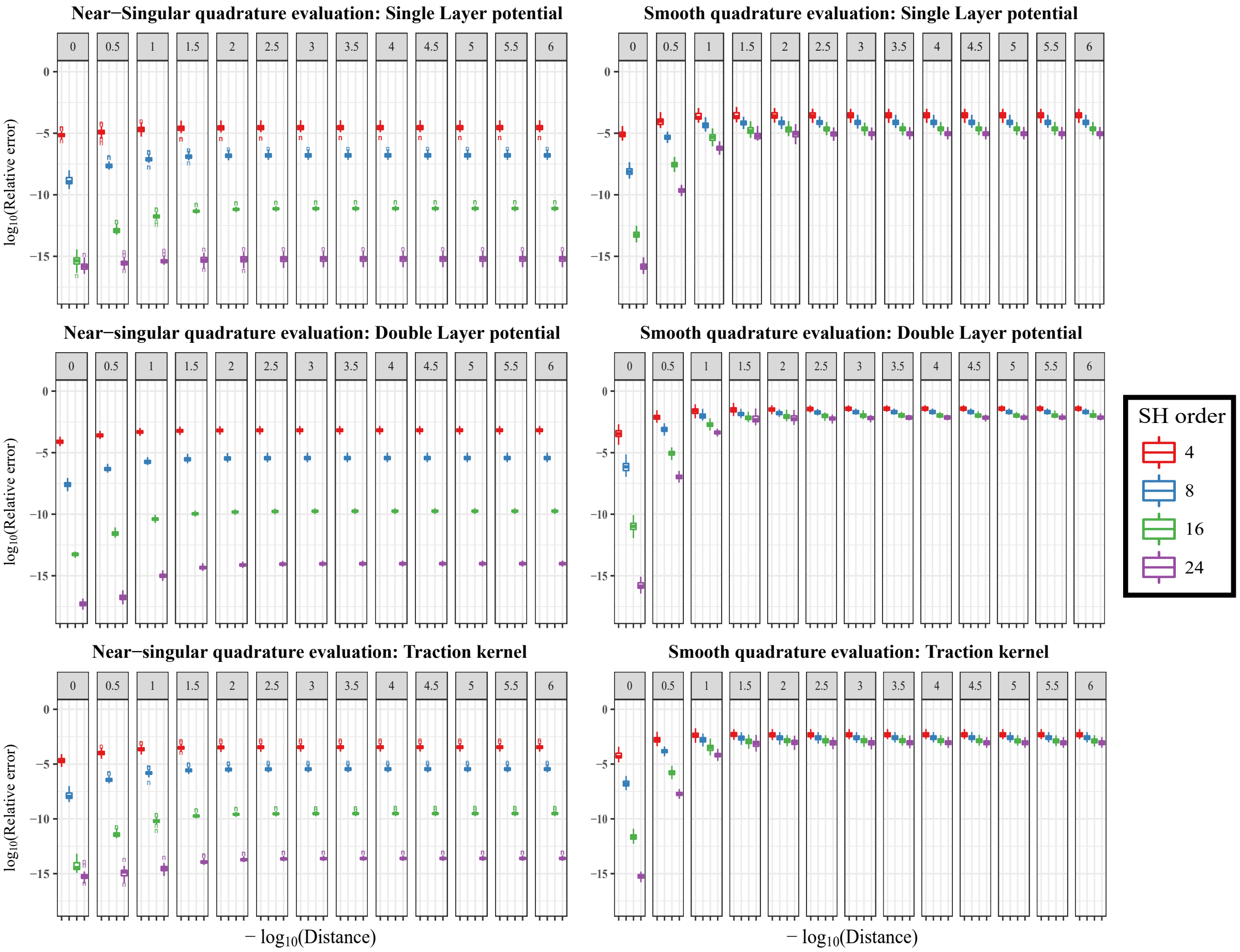

Figure 7: Quadrature performance on BIE data. *We measure the* log *relative error for the evaluation of* $\mathcal{S}$,$\mathcal{D}$ *and* $\mathcal{K}$ *for the Stokes equation for both integration schemes. For spherical harmonic order* $p = 4, 8, 16, 24$*, we present box plots categorized by* $-\log$ *of the distance to the sphere surface divided by its radius.*

We present the resulting box plot in Fig. 7. We again observe spectral convergence and robust accuracy for the near-singular scheme up to the boundary, requiring an expansion order of $p = 24$ to reach near machine-precision. In this and related experiments, we observe that as we bring the spheres closer, or bring the singular point sources closer to the sphere boundaries, the order required to achieve a desired target accuracy grows. We also observe in the experiments conducted, evaluation of $\mathcal{S}$ is typically a bit more accurate than that of $\mathcal{D}$ and $\mathcal{K}$.

### Scaling tests

Finally, we wish to test how the computational cost of our fast evaluation algorithm, as reflected by matrix apply times, scales as we increase the number of spheres $n_b$ for a given boundary integral operator. As mentioned in §4, FMM 3D libraries are employed to evaluate far-interactions using the smooth quadrature. All experiments are run serially on the Flux HPC cluster at University of Michigan.

For this purpose, we first produce two sets of monodisperse (all radii are equal), cubic lattices with spheres on each vertex of radii $r_v = (1 - 2^{-q})$, for $q = 1$ and $q = 4$. The volume fraction these occupy is $\frac{\pi}{6}(1 - 2^{-q})^3$, and minimum distance is $2^{-q+1}$. In order to test polydisperse suspensions, we add one sphere on each face of the cube, of radii $r_f = (1 - 2^{-q})(2 - \sqrt{2})$, and one at the center with $r_c = (1 - 2^{-q})(\sqrt{2} - 1)$. The volume fraction is increased to $\frac{\pi}{6}(1 - 2^{-q})^3 \left(1 + 2(2 - \sqrt{2})^3 + (\sqrt{2} - 1)^3\right)$, and minimum distance decreases to. We present the basic units in Fig. 8.

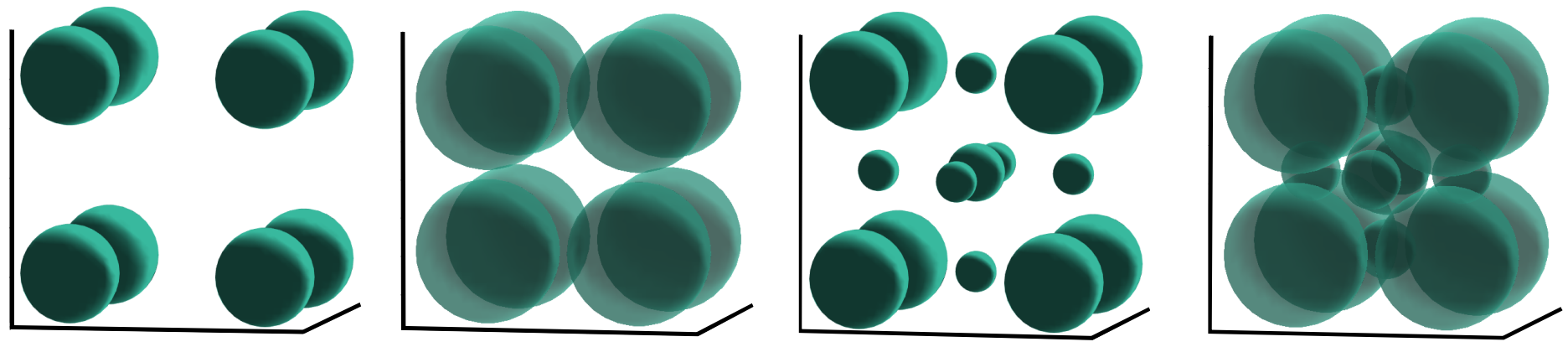

(a) monodisperse, q=1 (b) monodisperse, q=4 (c) polydisperse, q=1 (d) polydisperse, q=4

Figure 8: Lattice units. *We present the main units in the regular lattices used in our scaling tests. Spheres at cube vertices are all of radii* $r_v = (1 - 2^{-q})$*, face spheres* $r_f = (1 - 2^{-q})(2 - \sqrt{2})$*, and center spheres* $r_c = (1 - 2^{-q})(\sqrt{2} - 1)$.

For lattices from $2 \times 2 \times 2$ to $16 \times 16 \times 16$ units, we record the average time it takes to evaluate the Stokes double layer potential and the Traction kernel for the single layer potential. We omit the evaluation times for the single layer kernel as they are very similar to those for the double layer. For monodisperse lattices, the resulting number of bodies, $n_b \in \{8, 27, 125, 729, 4913\}$. For the polydisperse lattices, $n_b \in \{15, 71, 429, 2969, 22065\}$. Recall that, for spherical harmonic order $p$, the total number of unknowns is $N = 6p(p + 1)n_b$.

In Fig. 9, we compare the log-log plots for the Matrix vector apply times as a function of the number of spheres $n_b$. Across the experiments performed, experimental scaling matches our expectations, with slopes across experiments approaching 1 as $n_b$ increases.

Comparing results, we can readily observe that evaluation of the traction kernel $\mathcal{K}$ is slightly more expensive than that of $\mathcal{D}$ or $\mathcal{S}$ due to the additional terms needed for the near-singular evaluation scheme, as shown in §3.3. For the loosely-packed lattices, since each sphere has few neighbors, most of the apply time and resulting scaling is due to the far-field computation performed by the FMM, and scaling is closer to $O(n_b p^2)$ when comparing results for $p = 4$ and $p = 8$. In contrast, self and near-field evaluation becomes significantly more expensive for the packed lattice experiments, resulting in an increase in apply times (as compared to the loosely packed case) as well as scaling closer to $O(n_b p^3)$.

# 7 Conclusions

We presented formulas for the spectra of various boundary integral operators, which are defined on the unit sphere and whose kernels are the fundamental solutions to the Stokes equations. These allow rapid application of BIOs (e.g., $\mathscr{O}(p^3 \log p)$ as opposed to $\mathscr{O}(p^4)$ cost typically required for arbitrary particle shapes [17]) and analysis of more general integro-differential operators that arise in particulate fluid mechanics (e.g., see [41] for vesicle flows). In addition, we presented analytical expressions for evaluating the BIOs away from the boundary of the unit sphere. These permit accurate computation of hydrodynamic interaction of particles even in the case of closely packed suspensions. We discussed how the formulas can be used for solving various problems of interest and presented results verifying

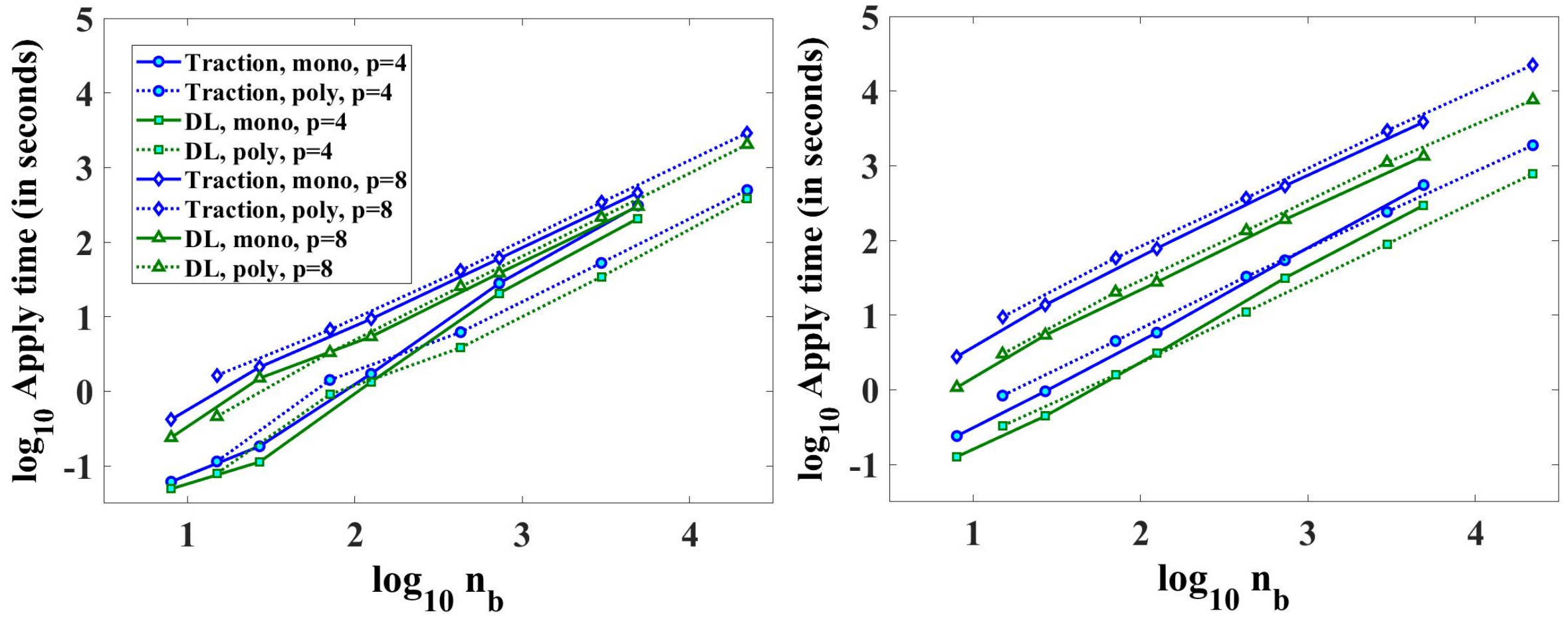

**Loosely packed lattices (q=1)** **Densely packed lattices (q=4)**

Figure 9: OPERATOR EVALUATION SCALING TESTS. *We measure the* log *apply time (in seconds) for the evaluation of* $\mathcal{D}$ *(in green) and* $\mathcal{K}$ *(in blue) for the Stokes equation for our integration scheme, for spherical harmonic order* $p = 4, 8$*. We plot it against* log *of the number of bodies* $n_b$ *for the loosely packed lattices* $(q = 1)$ *on the left, and for the densely packed lattices* $(q = 4)$ *on the right.*

the accuracy and scaling of our algorithm. These formulas may also be useful in developing high-order perturbation solution techniques for deformable particle flows similar to the ideas discussed in Vlahovska et al. [45], which we plan to explore in the future.

One important issue we haven't addressed in this work is that when particles are in close-to-touching regime, the density of the layer potentials becomes singular. This, in turn, means that the number of terms needed to represent the density function increases rapidly. Two promising approaches that ameliorate this problem are the hybrid method of images approach [5, 14] and the so-called *Recursively Compressed Inverse Preconditioning* (RCIP) method [20]. However, they have not been applied to the solution of the Stokes equations and moreover RCIP was applied for problems in 2D only; we plan to investigate their scope for the problems considered here.

# 8 Acknowledgements

We acknowledge support from NSF under grants DMS-1454010 and DMS-1418964. This research was supported in part through computational resources and services provided by Advanced Research Computing (ARC) at the University of Michigan, Ann Arbor. Part of this work was done while SV was visiting the Flatiron Institute; he is grateful for the hospitality and support received from the Biophysical Modeling group.

# A Formula derivations of Stokes potentials

We present a general derivation of the formulas in Eqs. 3.6 – 3.8 and 3.9 – 3.11 for the Stokes single and double layer potentials. We note that a version of this focused on the single layer potential can be found in [42].

As mentioned in §3.2, for given $n, m$, an ansatz velocity field of the form $\boldsymbol{u} = f(r)\boldsymbol{V}_n^m + g(r)\boldsymbol{W}_n^m + h(r)\boldsymbol{X}_n^m$ with associated pressures of the form $p_n^m = \frac{q(r)}{r} Y_n^m$ is proposed, and we require that it satisfies the Stokes equations. We arrive at the following system of four homogeneous ODEs by substituting in

the Stokes equations:

$$r^2 f_{rr} + 2r f_r - (n+1)(n+2)f + \frac{1}{2n+1}(-rq_r + (n+1)q) = 0 \tag{A.1}$$

$$r^2 g_{rr} + 2r g_r - n(n-1)g + \frac{1}{2n+1}(rq_r + nq) = 0 \tag{A.2}$$

$$r^2 h_{rr} + 2r h_r - n(n+1)h = 0 \tag{A.3}$$

$$(n+1)r f_r + (n+1)(n+2)f - nr g_r + n(n-1)g = 0 \tag{A.4}$$

of which the first three correspond to the momentum equation, and the last one to the continuity equation. The analytical solution of this system produces six sets of solutions,

| | $(i)$ | $(ii)$ | $(iii)$ | $(iv)$ | $(v)$ | $(vi)$ |
|---|---|---|---|---|---|---|
| $f_n(r)$ | $r^{-n-2}$ | $\frac{-n}{n+1}r^{-n}$ | 0 | 0 | $\frac{2n}{(2n+1)(2n+3)(n+1)}r^{n+1}$ | 0 |
| $g_n(r)$ | 0 | $\frac{2}{2n-1}r^{-n}$ | 0 | $r^{n-1}$ | $\frac{1}{2n+1}r^{n+1}$ | 0 |
| $h_n(r)$ | 0 | 0 | $r^{-n-1}$ | 0 | 0 | $r^n$ |
| $q_n(r)$ | 0 | $\frac{-2n(2n+1)}{n+1}r^{-n}$ | 0 | 0 | $-2r^{n+1}$ | 0 |

(A.5)

from which $(i), (ii)$ and $(iii)$ are only admissible in the exterior $(r > 1)$ (the rest are unbounded), and $(iv), (v)$ and $(vi)$ are only admissible in the interior $(r < 1)$, as the other three solutions are singular at the origin. Thus, a solution for the velocity field $\boldsymbol{u}$ and pressure $p$ may be written as a linear combination of three solutions in the exterior and three in the interior of the unit sphere. In order to determine the six resulting unknowns, we enforce the corresponding jump conditions for the given layer potential $\boldsymbol{u}$ and its traction $T[u] = (\nabla \boldsymbol{u} + \nabla \boldsymbol{u}^T)\boldsymbol{n} - p\boldsymbol{n}$ at the sphere,

$$[[\mathcal{S}[\boldsymbol{\sigma}]]]_\Gamma = 0 \quad \text{and} \quad [[T[\mathcal{S}[\boldsymbol{\sigma}]]]]_\Gamma = -\boldsymbol{\sigma}. \tag{A.6}$$

$$[[\mathcal{D}[\boldsymbol{\sigma}]]]_\Gamma = \sigma \quad \text{and} \quad [[T[\mathcal{D}[\boldsymbol{\sigma}]]]]_\Gamma = 0. \tag{A.7}$$

Substituting this linear combination and taking the corresponding inner products with $\boldsymbol{V}_n^m$, $\boldsymbol{W}_n^m$ and $\boldsymbol{X}_n^m$ yields three $6 \times 6$ linear systems for each layer potential. We note that since the equations in A.4 for $f$ and $g$ aren't coupled with that for $h$, these reduce to one $4 \times 4$ ($(i), (ii), (iv)$ and $(v)$) and one $2 \times 2$ ($(iii), (vi)$ independent linear systems. Solving these for arbitrary $n$ provides the desired formulas for both layer potentials. We note that in this calculation, it's useful to employ the formula in Eq. (3.18), which for a velocity field $\boldsymbol{u} = f(r)\boldsymbol{V}_n^m + g(r)\boldsymbol{W}_n^m + h(r)\boldsymbol{X}_n^m$ and pressure $p = \frac{q(r)}{r}Y_n^m$ simplifies on the sphere to:

$$\begin{aligned}&\left(\frac{(3n+2)f_r(1) - n(n+2)f(1) - ng_r(1) + n(n-1)g(1) - q(1)}{2n+1}\right)\boldsymbol{V}_n^m \\ &+ \left(\frac{-(n+1)f_r(1) - (n+1)(n+2)f(1) + (3n+1)g_r(1) + (n+1)(n-1)g(1) + q(1)}{2n+1}\right)\boldsymbol{W}_n^m \\ &+ (h_r(1) + h(1))\,\boldsymbol{X}_n^m.\end{aligned} \tag{A.8}$$